\documentclass[12pt,eqno]{article}
\usepackage{amsmath}
\usepackage{amsfonts}
\usepackage{mathrsfs}
\usepackage{amssymb,amsmath,latexsym}

 \usepackage[dvips]{color}
\usepackage{amsmath,amsthm}     
\usepackage{amssymb}            
\usepackage{euscript}           
\usepackage{graphicx}
\usepackage{tikz}
\usepackage{curves,enumerate}

\usepackage[colorlinks]{hyperref}

\newtheorem{theorem}{Theorem}[section]

\newtheorem{lemma}[theorem]{Lemma}
\newtheorem{rem}{Remark}
\newtheorem{example}{Example}[section]
\newtheorem{corollary}[theorem]{Corollary}
\newtheorem{definition}{Definition}[section]

\newtheorem{hyp}{Hypothesis}

\def\lab#1{\mbox{$\label{#1}$}}

\def\ga{\gamma}                         

\def\ga{\gamma}
\def\ra{\rightarrow}                    
\def\Ra{\Rightarrow}                    
\def\proof{\noindent {\bf Proof: }\ }

\def\wt{\widetilde}
\def\qed{\hfill $\square$ \bigskip}

\def\lab{\label}
\def\beq{\begin{equation}}               
\def\eeq{\end{equation}}                 
\def\bea{\begin{eqnarray}}             
\def\eea{\end{eqnarray}}               
\def\be*{\begin{eqnarray*}}             
\def\ee*{\end{eqnarray*}}               
\def\ba{\begin{array}}                  
\def\ea{\end{array}}                    

\def\;{\vspace{3mm} \\ }
\def\Q{\rm Q}
\def\N{\mathbb{N} }

\def\R{\mathbb{R}}
\def\E{\mathbb{E} }
\def\P{\mathbb{P} }
\def\r{\mathbb{R}}

\def\~{\widetilde}

\def\beqlb{\begin{eqnarray}} \def\eeqlb{\end{eqnarray}}
\def\beqnn{\begin{eqnarray*}} \def\eeqnn{\end{eqnarray*}}

\def\<{\langle}  \def\>{\rangle}

\def\wt{\widetilde}

\def\bB{\mbox{\boldmath $B$}} 
 \def\B{{\cal B}}

\def\ra{\rightarrow}                    
\def\Ra{\Rightarrow}                    

\def\bde{\begin{definition}}
\def\ede{\end{definition}}

\def\bth{\begin{theorem}}
\def\eth{\end{theorem}}
\def\ble{\begin{lemma}}
\def\ele{\end{lemma}}
\def\bcor{\begin{corollary}}
\def\ecor{\end{corollary}}
\def\bhyp{\begin{hyp}}
\def\ehyp{\end{hyp}}
\def\brem{\begin{rem}}
\def\erem{\end{rem}}

\begin{document}


\title{\Large {\bf Tanaka formula and local time for a class of interacting branching measure-valued diffusions}
\thanks{Partial funding in support of this work was provided by the Natural
Sciences and Engineering Research Council of Canada (NSERC) and the Department of Mathematics at the University of Oregon.}}

\author{\bf D. A. Dawson\thanks{D. A. Dawson, Carleton University, ddawson@math.carleton.ca}
\and \bf J. Vaillancourt\thanks{J. Vaillancourt, HEC Montr\'eal, jean.vaillancourt@hec.ca.
Corresponding Author : 3000, chemin de la C\^ote-Sainte-Catherine, Montr\'eal (Qu\'ebec),
                                                      Canada H3T 2A7 orcid.org/0000-0002-9236-7728 }
\and \bf H. Wang\thanks{H. Wang, University of Oregon, haowang@uoregon.edu}}

\maketitle


{\narrower{\narrower

\centerline{\bf Abstract}

\medskip

We construct superprocesses with dependent spatial motion (SDSMs) in Euclidean spaces $\R^d$ with $d\ge1$ and show that,
even when they start at some unbounded initial positive Radon measure such as Lebesgue measure on $\R^d$, their
local times exist when $d\le3$. A Tanaka formula of the local time is also derived.
\medskip



\par}\par}

 \baselineskip=6.0mm

 \section{Introduction}

 \setcounter{equation}{0}

In the present paper, we provide an explicit representation of the local time for a class of interacting
superprocesses (measure-valued branching diffusions) on Euclidean space $\R^d$ when $d\le3$.
Known as superprocesses with dependent spatial motion (SDSM), this class was introduced in
Wang \cite{Wang97}, \cite{Wang98} and extended in Dawson et al. \cite{DLW01}. 

The subdivision of SDSM into subclasses, according to the portion of the state space that SDSM charges almost surely, its
various representations by associated stochastic partial differential equations (SPDE)
and some of its trajectorial properties have been exhibited and analyzed since then,
notably in the special case of measure-valued processes on the real line ($d=1$).
In that case, let us mention the following results:
the dimension of its support in Wang \cite{Wang97}, depending on the degeneracy or not of the diffusion term;
the explicit form of the SPDE for the motion
of the resulting purely atomic measure-valued SDSM in various degenerate cases,
in Dawson et al. \cite{DLW03} and Li et al. \cite{LWX04a};
an explicit representation of its density in Dawson et al. \cite{DVW2000} in the non degenerate case
and the delicate matter of joint continuity in time and space of this density,
in Li et al. \cite{LWXZ12} and Hu et al. \cite{HLN13} and \cite{HNX19}. 
Extending some of these results to the general case $d\ge1$ still offers many challenges, some of which are addressed here. 

To avoid repetitions, we make the following basic assumptions, called for as needed throughout this paper and
related to the properties of the processes themselves, as well as the filtered probability spaces they are constructed on. 
Definitions and notations are gathered at the beginning of Section \ref{sec:main}; they are introduced in the text by the $:=$ symbol. 

\begin{hyp}\label{hyp:basicassumpFilter}
Let $(\Omega, {\cal F}, \{{\cal F}_t\}_{t \geq 0}, \P)$ be a filtered probability space with a right continuous filtration
$\{{\cal F}_t\}_{t \geq 0}$, satisfying the usual hypotheses and upon which all our processes are built, notably an
$\R^1$-valued Brownian sheet $W$ on $\R^d$ (see below) and a countable family $\{B_{k}, k\geq 1\}$ of independent,
$\R^d$-valued, standard Brownian motions written $B_{k}=(B_{k1}, \cdots, B_{kd})$.
The family $\{B_{k}, k\geq 1\}$ is assumed independent of $W$.
\end{hyp}
Following Walsh \cite[Chapter 2]{Walsh86}, a random set function $W$ on ${\cal B}(\R^d \times[0,\infty))$ defined on
$(\Omega, {\cal F}, \{{\cal F}_t\}_{t \geq 0}, \P)$ is called an $\R^1$-valued Brownian sheet on $\R^d$
(or space-time white noise) if both of the following statements hold:
for every $A \in {\cal B}(\R^d)$ having finite Lebesgue measure $\lambda_0(A)$,
the process $M(A)_t:=W(A\times [0, t])$ is a square-integrable $\{ {\cal F}_t\}$-martingale;
and for every pair $A_i\in {\cal B}(\R^d \times[0,\infty))$, $i=1, 2$, having finite
Lebesgue measure with $A_1\cap A_2=\emptyset$, the random variables $W(A_1)$ and
$W(A_2)$ are independent, Gaussian random variables with mean zero,
respective variance $\lambda_0(A_i)$ and
$ W(A_1\cup A_2) = W(A_1)+W(A_2)$ holds $\P$-almost surely
(see Walsh \cite{Walsh86}, Dawson \cite[Section 7.1]{Dawson93}
 and Perkins \cite{Perkins02} for further details). 

The rest of our assumptions relate to the coefficients in the equations (bounded and Lipschitz continuous) 
and the restrictions imposed on the (often infinite) initial measure.

 \begin{hyp}\label{hyp:basicassumpElliptic} 
The vector $h=(h_1, \cdots,h_d)$ satisfies $h_p \in  L^1(\R^d) \cap {\rm Lip}_b  (\R^d)$
and the $d\times d$ matrix $c=(c_{pr})$ satisfies $c_{pr} \in   {\rm Lip}_b(\r^d)$, for $p,r=1, \cdots,d$. 
For each $p,q=1,\ldots,d$, we write $a_{pq}(x):=\sum_{r=1}^{d}c_{pr}(x)c_{qr}(x)$ and 
$\rho_{pq}(x-y):=\int_{\R^d}h_{p}(u-x)h_{q}(u-y)du$.  
 For every $m\ge1$, the $dm\times dm$ diffusion matrices $(\Gamma_{pq}^{ij})_{1 \leq i,j \leq m; 1 \leq p,q \leq d}$
 of real-valued functions defined by
  \bea \label{gammaij} \hspace*{8mm}\Gamma_{pq}^{ij}(x_{1}, \cdots,
x_{m}) := \left\{
\begin{array}{ll}
                                 (a_{pq}(x_{i})+\rho_{pq}(0))
                                 \quad & \mbox{if
                                 $i=j$,} \\
                                 \rho_{pq}(x_{i}- x_{j}) \quad & \mbox {if $i\neq j$},
                                 \end{array}
                                 \right.
\eea
are strictly positive definite everywhere on $\R^d$ with
smallest and largest eigenvalues bounded away respectively from $0$ and $\infty$, uniformly in $\R^d$;
that is, there are two positive constants $\lambda_m^*$ and $\Lambda_m^*$ such that for any
$\xi = (\xi^{(1)}, \cdots, \xi^{(m)}) \in (\R^d)^m$ we have a positive definite form which satisfies
\[
 0 <  \lambda_m^* |\xi|^2\leq \sum_{k,l=1}^m \sum_{p,q=1}^d \Gamma_{pq}^{kl}(\cdot) \xi_p^{(k)} \xi_q^{(l)}
 \leq \Lambda_m^*|\xi|^2 < \infty.
\]
\end{hyp} 

 \noindent  {\bf Remark:} 
Hypothesis \ref{hyp:basicassumpElliptic} ensures that $\Gamma$ is uniformly elliptic (lower bound). 
Since $\Gamma$ is a sum of two non-negative definite matrices, this occurs as soon as one of the two 
matrices is strictly positive definite. This is the case if the individual motions are uniformly elliptic. 
Consult Section \ref{sec:SDSM} for more consequences of Hypothesis \ref{hyp:basicassumpElliptic}. 

Hypotheses \ref{hyp:basicassumpFilter} and \ref{hyp:basicassumpElliptic} together  
guarantee the existence of non-degenerate SDSM and their characterization 
by way of the martingale problem approach (Theorem \ref{MPforSDSM}). 

The rest of our assumptions relate to the restrictions imposed on the initial measure. 
Their statement requires some additionnal notation. 

Let $M(\R^d)$ be the set of all positive Radon measures on $\R^d$  and
$M_0(\R^d)$, its subspace of finite positive Radon measures.
For any $a \ge 0$, let $I_a(x):= (1 + |x|^2)^{(-a/2)}$ and 
\bea \label{spacesMa}
M_a(\R^d)= \{\mu \in M(\R^d): \<I_a,\mu\>:= \int_{\R^d}I_a(x) \mu(dx) < \infty \}.
\eea
The topology $\tau_a$ of $M_a(\R^d)$ is defined in the following way: $\mu_n \in M_a(\R^d)$ converges to $\mu \in M_a(\R^d)$ 
as $n \rightarrow \infty$, if $\lim_{n \rightarrow \infty}\< \phi, \mu_n\> = \<\phi, \mu\>$  holds for every $\phi \in K_a(\R^d),$ where
\bea \label{spacesKa}
K_a(\R^d)=\{\phi: \phi:= h + \beta I_a, \beta \in \R, h \in C_c(\R^d)\}.
\eea
Then, $(M_a(\R^d), \tau_a)$ is a Polish space
(see Iscoe \cite{Iscoe86a} and Konno and Shiga \cite{KonnoShiga88}).
For instance, the Lebesgue measure $\lambda_0$ on $\R^d$ belongs to $M_a(\R^d)$ for any $a >d$.
Furthermore, both $dx=\lambda_0(dx)$ are used interchangeably when calculating Lebesgue integrals. 

\begin{hyp}\label{hyp:basicassumpGauss}
For any $T > 0$, the initial measure $\mu_0\in\cup_{a\ge0}M_a(\R^d)$ verifies
\bea \label{nonBA}
 \sup_{x\in \R^d}  \sup_{0 < t \leq T} \<\varphi_t(x-\cdot), \mu_0\>  < \infty,
\eea
where $\varphi_t$ is the transition density of the standard $d$-dimensional Brownian motion.
\end{hyp}

In order to fulfill the requirements of using Fubini's theorem in some of the proofs, we also require an additional uniform bound
on measure $\mu_0 \in M_a(\R^d)$. 

\begin{hyp}\label{hyp:basicassumpUniformInteg} 
The initial measure $\mu_0 \in M_a(\R^d)$ satisfies 
\[
 \sup_{w\in \R^d}  \<I_a(\cdot+w), \mu_0 \> < \infty .
\] 
\end{hyp} 

 \noindent  {\bf Remark:} 
Hypotheses \ref{hyp:basicassumpGauss} and \ref{hyp:basicassumpUniformInteg} are uniform boundedness conditions ensuring the existence of the local time for SDSM in general (Theorem \ref{lt_th1}) as well as the validity of the Tanaka formula (\ref{TanakaII}).  
Hypothesis \ref{hyp:basicassumpGauss} is akin to condition (2.9) in Sugitani \cite{Sugitani89}, instrumental in his proof of joint continuity of the local time for Super-Brownian motion. The treatment of the joint continuity of the local time for SDSM in general requires sharper inequalities than the ones used here and is treated in \cite{DVW2021}. 
Additional insight into the need for Hypotheses \ref{hyp:basicassumpGauss} and \ref{hyp:basicassumpUniformInteg} 
is provided in Section \ref{sec:SDSM}. 

The class SDSM $\{\mu_t\}$ is the main object of this paper and is characterized through the second order linear differential operator
${\cal L}$, defined for bounded smooth real-valued functions $F$ on $M_a(\R^d)$ by  

   \beqlb \lab{pregenerator}
    {\cal L}F(\mu) & := & { \displaystyle \frac{1}{2} \sum_{p,q=1}^{d} \int_{\R^d}  (a_{pq}(x) + \rho_{pq}(0)) 
  \left(\frac{\partial^{2}}{\partial x_p \partial x_q}\right) 
  \frac{{\delta} F(\mu)}{{\delta} \mu(x)}\,\mu(dx)} \nonumber \\
  & & + { \displaystyle \frac{1}{2}\sum_{p,q=1}^{d} \int_{{\R^d}}  \int_ {{\R^d}}
 \rho_{pq}(x-y) \left(\frac{\partial}{\partial x_p }\right)
 \left(\frac{\partial}{\partial y_q }\right)
 \frac{{\delta}^{2}F(\mu)}{{\delta} \mu(x) {\delta} \mu(y)} \mu(dx) \mu(dy)} \nonumber \\ 
  & & + { \displaystyle \frac{\gamma \sigma^2}{2}
 \int_{{\R^d}} \frac{{\delta}^{2} F(\mu)}{{\delta} \mu(x)^{2}} \mu (dx). }
   \eeqlb

 The local (or individual) diffusion coefficients $a_{pq}$ and the global (or common) interactive diffusion coefficients $\rho_{pq}$ 
 are defined within Hypothesis \ref{hyp:basicassumpElliptic}. 
 Parameter $\gamma>0$ is related to the branching rate of the particle system
 and $\sigma^2 > 0$ is the variance of  the limiting offspring distribution. The variational derivative is defined
  by
  \[ \frac{\delta F(\mu)}{\delta \mu(x) }:= \lim_{\epsilon
 \downarrow 0} \frac{F(\mu+\epsilon \delta_{x})-F(\mu)}{\epsilon}
  \]
 where $\delta_{x}$ stands for the Dirac measure at $x$. The domain ${\cal D}({\cal L})$ of the operator ${\cal L}$
 contains all functions of the form $F(\mu) = g(\<\phi_1, \mu\>,  \cdots, \<\phi_k, \mu \>)$
with $g \in C_b^2(\R^k)$ for some $k\ge1$ and $\phi_i  \in C_c^{\infty}(\R^d)$ for every $1 \leq i \leq k $ (although it is not 
restricted to just these functions). For any $\mu \in M(\R^d)$ and any $\mu$-integrable $\phi$ 
we write $\<\phi,\mu\>= \int_{\R^d} \phi(x) \mu(dx)$ here and henceforth.

Clearly the class of SDSM includes the critical branching Dawson-Watanabe superprocesses when $h\equiv0$.
The literature on these is extensive and the reader may consult the lecture notes by Dawson \cite{Dawson92, Dawson93},
Le Gall \cite{LeGall99}, Etheridge \cite{Etheridge00} and Perkins \cite{Perkins02} for historical insights into the origin 
and the early evolution of the field, 
as well as the subsequent works by Li \cite{Li10} and Xiong \cite{Xiong13} for thorough updates on the subject.

Amongst the many properties of SDSM the one of interest here is the existence of a local time when $d\le3$. 
A local time of SDSM $\{\mu_t\}$ is a density process of the occupation time
 process $\int_0^t \mu_s ds$ of SDSM, a time-averaging giving rise to a new measure-valued process
 with more regular paths and, in some cases, a density with respect to Lebesgue measure, even when SDSM
 itself does not have one. For instance, from Wang \cite{Wang97} we know that in the degenerate case, the
 SDSM is a purely atomic measure-valued process, so the density of the occupation time for this degenerate SDSM process
 may not exist. However, Li and Xiong \cite{LiXiong07} introduced an interesting alternative way to define the local time for a class of
 purely atomic measure-valued processes along the path of each particle. The local time (in this sense) of the degenerate SDSM
 is constructed there and its joint H{\"o}lder continuity proved. We now proceed in our non-degenerate (uniformly elliptic) case with the following,
 more familiar definition.
 \bde \lab{de1}
 A local time of the SDSM $\{\mu_t\}$ is any ${\cal B}(\R^d \times[0,\infty)) \times {\cal F}$-measurable function 
  $\Lambda^{x}_t(w):(\R^d \times [0, \infty) \times \Omega) \rightarrow [0, \infty)$, 
  satisfying 
 \beqlb
 \label{LT}
 \int_{\R^d}\phi(x) \Lambda_t^x dx = \int_0^t \<\phi, \mu_s\>ds,  
 \eeqlb 
 for every $t > 0$ and every $\phi \in C_c(\R^d)$. 
\ede 
We will construct such a local time so that $(\ref{LT})$ holds simultaneously for every $t > 0$ 
and every $\phi \in C_c(\R^d)$ outside of a common $\P$-null set. 

In the degenerate case, Li et al. \cite{LWX04a} have shown that two SDSM with the same initial data can have distinct pathwise solutions, 
thus generating two distinct local times that satisfy equation $(\ref{LT})$; hence the need for the alternate definition of local time provided in 
Li and Xiong \cite{LiXiong07} in the degenerate case. 
Here instead, under Hypotheses \ref{hyp:basicassumpFilter} through \ref{hyp:basicassumpUniformInteg}, 
this situation is avoided, as will be seen in Section \ref{sec:main}. 

In the case of Super-Brownian motion (where $h\equiv0$ and $c$ is the identity matrix) the existence and the joint space-time continuity of paths for its local time when $d\le3$ go back to Iscoe \cite{Iscoe86b} and Sugitani \cite{Sugitani89}. These results,
as well as further path properties, were generalized to superdiffusions (still $h\equiv0$) in Krone \cite{Krone93}. In these and many other papers where the finer aspects of the superprocesses are analyzed, the argument largely depends on a multiplicative property of branching processes and the availability of a manageable closed form for the log-Laplace functional, a powerful tool to estimate the higher moments of $\{\mu_t\}$. Unfortunately, in our model the dependency of motion ($h$ is no longer the null function) destroys the multiplicative property in question and makes this approach largely intractable, as it relies intimately on the independence structure built into Dawson-Watanabe superprocesses. This method was applied notably by Adler and Lewin \cite{AdlerLewin92} in their proof of the Tanaka formula for the local time of Super-Brownian motion and super stable processes. This also occurs when trying the approach proposed in
L\'opez-Mimbela and Villa \cite{Lopez-MimbelaVilla04} for Super-Brownian motion, where an alternative representation of the local time simplifies the proof of its joint continuity by taking advantage of sharp estimates for the Green function of Brownian motion and its associated Tanaka formula. However, the higher order singularity of the Green function and its derivative in our case,
raises some new technical difficulties in the moment estimation of the interacting term, as well as in the handling of a stochastic convolution integral term appearing in the corresponding Tanaka formula. These issues are resolved using a perturbation argument and a Tanaka formula for SDSM emerges. 

The remainder of this paper is organized as follows. 
Section \ref{sec:main} states the main results and assembles the remaining notation required for their formulation. 
Section \ref{sec:SDSM} is devoted to the construction of SDSM $\{\mu_t\}$ started with an unbounded initial measure 
on $\R^d$ with $d \geq 1$. This construction proceeds as a limit of a sequence of branching particle systems in $\R^d$; 
it uses a tightness argument adapted from the one built by 
Ren et al. \cite{RenSongWang09} for SDSM on a bounded domain $D \subset \R^d$. 
The derivation of the associated SPDE is achieved there as well. 
In Section \ref{sec:dualConst}, sharp bounds are obtained for the $k^{th}$-moments of SDSM 
by computing them through a duality argument, using the martingale problem approach.
In Section \ref{sec:Tanaka}, we derive a Tanaka formula for SDSM and
use it to prove the existence of the local time $\Lambda_t^x$ for SDSM for some unbounded initial measures. 
To avoid forward referencing, Sections \ref{sec:SDSM}, \ref{sec:dualConst} and \ref{sec:Tanaka} are in logical and chronological order.
Some technical results have their proof postponed to Section \ref{app:ProofsOfLemmas}; those proofs are self-contained. 

\section{Main results and notation} \label{sec:main}
\setcounter{equation}{0} 

For any Polish space $S$, that is, a topologically complete and separable metric space,
${\cal B}(S)$ denotes its Borel $\sigma$-field, $B(S)$ the Banach space
of real-valued bounded Borel measurable functions on $S$ with the supremum norm $\|\cdot\|_{\infty}$
and $C(S)$ the space of real-valued continuous functions on $S$.
Subscripts $b$ or $c$ on any space of functions will always refer to its subspace of bounded
or compactly supported functions, respectively, as in $C_b(S)$ and $C_c(S)$ here.
$S^m$ denotes the $m$-fold product of $S$.

The spaces of continuous $C([0, \infty), S)$ and c\`adl\`ag $D([0, \infty), S)$
trajectories into Polish space $S$ are respectively equipped with the topology of uniform convergence on compact time sets
and the usual Skorohod topology; they are themselves also Polish spaces (see Ethier and Kurtz \cite{EthierKurtz86}).

Given any positive Radon measure $\mu\in M(\R^d)$ and any $p\in[1,\infty)$,
 we write $L^p(\mu)$ for the Banach space of real-valued Borel measurable functions on $\R^d$,
 with finite norm $\|\phi\|_{\mu,p}:= \{\int_{\R^d}|\phi(x)|^pd\mu(x)\}^{1/p} < \infty$
 and $|x|^2 := \sum_{i=1}^dx_i^2$. When $\mu=\lambda_0$ is the Lebesgue measure
 we use the standard notation $L^p(\R^d)=L^p(\lambda_0)$ and $\|\phi\|_p:= \|\phi\|_{\lambda_0,p}$.

 We need various subspaces of continuous functions inside $C(\R^d)$, notably
 $C^k(\R^d)$ the space of continuous functions on $\R^d$ with continuous derivatives up to and including order $k\ge0$,
with $C^\infty(\R^d)$ their common intersection (the smooth functions) and noticing that $C^0(\R^d)=C(\R^d)$;
$C_b^k(\R^d)$ their respective subspace of bounded continuous functions with bounded derivatives
up to and including order $k$, again with $C_b^\infty(\R^d)$ their common intersection and $C_b^0(\R^d)=C_b(\R^d)$;
$C_0^k(\R^d)$ those bounded continuous functions vanishing at $\infty$ together with their derivatives up to and including order $k$,
with $C_0^\infty(\R^d)$ their common intersection and
$C_0^0(\R^d)=C_0(\R^d)$, this last a Banach space when equipped with finite supremum norm;
$C_c^k(\R^d)$ the further subspace of those with compact support,
again with $C_c^\infty(\R^d)$ their common intersection and
$C_c^0(\R^d)=C_c(\R^d)$.
We use ${\rm Lip}(\R^d)$ to denote the space
of Lipschitz functions on $\R^d$, that is, $\phi \in {\rm Lip}(\R^d)$
if there is a constant $M>0$ such that $|\phi(x)-\phi(y)|\leq M |x-y|$
for every $x, y \in \R^d$. The class of bounded Lipschitz functions
on $\R^d$ will be denoted by ${\rm Lip}_b(\R^d)$. 
Let ${\cal S}(\R^d)$ be the Schwartz space of rapidly decreasing test functions and 
${\cal S}'(\R^d)$ the space of Schwartz tempered distributions, the dual space of ${\cal S}(\R^d)$  
(see Al-Gwaiz \cite{Al-Gwaiz92}, Barros-Neto \cite{Barros-Neto73} or Schwartz \cite{Schwartz59}).

In addition to ${\cal S}(\R^d)$, the main set of functions of interest here is $K_a(\R^d)$ defined in $(\ref{spacesKa})$ 
for any real number $a\ge0$. Since $C_c^\infty(\R^d)$ is uniformly dense in $C_c(\R^d)$ (with $C_0(\R^d)$ as common closure), 
the uniform closure of $K_a(\R^d)$ remains unchanged if we replace $C_c(\R^d)$ by $C_c^\infty(\R^d)$.
Both are also $\|\cdot\|_p$-dense in $L^p(\R^d)$ for every $p\in[1,\infty)$ 
(for instance, see Lemma 2.19 in Lieb and Loss \cite{LiebLoss01}), a fact that will come in handy later.  
Of course ${\cal S}(\R^d)$ is uniformly dense in $C_0(\R^d)$ and $\|\cdot\|_p$-dense in $L^p(\R^d)$ for every $p\in[1,\infty)$ as well.

We will also need $C_b^{1,2}([0, t] \times (\R^d)^m)$, the space of 
bounded continuous functions with all derivatives bounded, up to and including order $1$ in the time variable up to time $t$
and order $2$ in the $md$ space variables, including mixed derivatives of that order. 
When no ambiguity is present we also write the partial derivatives 
(of functions and distributions) in abridged form
\begin{equation*}
\partial_p:=\frac{\partial}{\partial x_p} \quad \mbox{and}\quad 
\partial_p\partial_q:=\frac{\partial}{\partial x_p}\frac{\partial}{\partial x_q} \quad \mbox{and so on.}
\end{equation*}

Set 
 \beqlb\label{eqn:Gn1}
G_1 := \sum_{p,q=1}^{d}\tfrac{1}{2}(a_{pq}(x)+\rho_{pq}(0)) \partial_p \partial_q.
 \eeqlb 

Let us now turn our attention to the characterization of SDSM through the formulation of a well-posed martingale problem 
(see Ethier and Kurtz \cite{EthierKurtz86}). 

A solution to the $({\cal L}, \delta_{{\mu}_{0}})$-martingale problem  
is a stochastic process $\mu$ with values in $M_a(\R^d)$ defined on 
$(\Omega, {\cal F}, \{{\cal F}_t\}_{t \geq 0}, \P)$ with initial value ${\mu}_{0}\in M_a(\R^d)$ such that, for every 
$F\in {\cal D}({\cal L})$, the process $F(\mu_t)-\int_{0}^{t} {\cal L} F(\mu_s)ds$ is an ${\cal F}_t$-martingale.
We say this martingale problem is well-posed (or has a unique solution) if such a solution exists and 
any two solutions have the same finite dimensional distributions.  


\bth \lab{MPforSDSM}
Assume Hypotheses $\ref{hyp:basicassumpFilter}$ and $\ref{hyp:basicassumpElliptic}$.
For any $a\ge0$ and initial value ${\mu}_{0}\in M_a(\R^d)$,
 the $({\cal L}, \delta_{{\mu}_{0}})$-martingale problem
 for the operator given by $(\ref{pregenerator})$ is well-posed and its
unique solution $\mu_t$ is an $M_a(\R^d)$-valued diffusion process which satisfies the SPDE
 \beqlb \label{SPDEc}
  \<\phi,\mu_{t}\>  -  \<\phi,\mu_{0}\>
  = {X}_{t}(\phi)   + M_t(\phi)  + \int_{0}^{t} \left\<  G_1\phi, \, {\mu}_{s} \right\> ds
 \eeqlb
$\P_{\mu_0}$-almost surely and the $\P_{\mu_0}$-null set in question is common to all $t>0$ 
and every $\phi \in K_a(\R^d) \cup {\cal S}(\R^d)$.  
 Further, for every $t>0$ and every $\phi \in K_a(\R^d) \cup {\cal S}(\R^d)$, there also holds 
  \beqlb \label{SPDEb}
   \E_{\mu_0} \bigg(\<\phi,\mu_{t}\>  -  \<\phi,\mu_{0}\>
  - {X}_{t}(\phi)   -  M_t(\phi)  - \int_{0}^{t} \left\<  G_1\phi, \, {\mu}_{s} \right\> ds\bigg)^2 = 0.
 \eeqlb
 Both
$$
{X}_{t}(\phi)  := \sum_{p=1}^{d} \int_{0}^{t}\int_{\r^d}
\left\<h_p(y-\cdot) \partial_p \phi(\cdot), {\mu}_{s} \right\>
{W}(dy,ds)
$$
and
$$ M_t(\phi) :=\int_0^t \int_{\R^d} \phi (y) M(dy, ds)
$$
are continuous square-integrable $\{{\cal F}_t\}$-martingales, mutually orthogonal for
every choice of $\phi \in K_a(\R^d) \cup {\cal S}(\R^d)$
and driven respectively by a Brownian sheet $W$ and a square-integrable martingale measure $M$ with
 \beqnn
  \<{M}(\phi) \>_t =
  \ga \sigma^2 \int_0^t \<\phi^2, {\mu}_s \> ds
  \qquad \hbox{for every } t>0 \hbox{ and } \phi \in K_a(\R^d) \cup {\cal S}(\R^d).
 \eeqnn
  \eth

 Here the filtration of choice is
 ${\cal F}_t:=\sigma\{ \<\phi,\mu_{s}\> , M_s(\phi), X_s(\phi): \phi \in K_a(\R^d), s \leq t\}$. 
The law for this solution $\{\mu_t\}$ started at ${\mu}_{0}\in M_a(\R^d)$ will henceforth be denoted by $\P_{\mu_0}$. 
The corresponding expectation is just $\E_{\mu_0}$. 

The proof of Theorem \ref{MPforSDSM} is in Section \ref{sec:dualConst}. 

This unique solution to the martingale problem for $({\cal L}, {\cal D}({\cal L}))$ is our SDSM. 

Note that Hypothesis \ref{hyp:basicassumpElliptic} ensures that operator $G_1$ is uniformly elliptic on $C^2(\R^d)$. 
For the single particle transition density $q_t^{1}(0,x)$ (from 0) of the semigroup $P_t^1$
associated with generator $G_1$ from (\ref{eqn:Gn1}),  
its Laplace transform (in the time variable) is given by
  \beqlb\label{eqn:green}
 Q^{\lambda}(x) := \int_0^{\infty}e^{-\lambda t} q_t^{1}(0,x)dt,
 \eeqlb
 for any $\lambda > 0$. Formally $Q^{0}$ is known as Green's function for density $q_t^{1}$
 and exhibits a potential singularity at $x=0$. For all $x \in \R^d \smallsetminus \{x=(x_1,x_2,\cdots,x_d):x_1x_2\cdots x_d=0\}$ 
 we can also write --- the legitimacy of the interchange and the finiteness of both expressions will be explained in due course ---
   \beqlb\label{eqn:interchange}
    \partial_{x_i}Q^{\lambda}(x)= 
    \partial_{x_i}\int_0^{\infty}e^{-\lambda t}q_t^{1}(0,x)dt
    =\int_0^{\infty}e^{-\lambda t}\partial_{x_i}q_t^{1}(0,x)dt 
 \eeqlb
  for any $i \in \{1,2,\cdots,d\}$, with the derivative taken in the distributional sense --- see Lemma \ref{FellerProp} for clarifications. 
The singularity at $x=0$ requires the following perturbation. For every $\epsilon>0$ and every $x \in \R^d$, define
 \[
  Q^{\lambda}_{\epsilon}(x - \cdot):=
   \int_{0}^{\infty} e^{- \lambda s}q_{s + \epsilon}^{1}(x- \cdot) ds
   =e^{\lambda \epsilon}\int_{\epsilon}^{\infty} e^{- \lambda t}q_t^{1}(x- \cdot) dt \in C_b^{\infty}(\R^d).
   \] 
   
The local time for SDSM is constructed by way of its approximating sequence, as follows: for every $(t,x) \in [0, \infty) \times \R^d$, define 
  \beqlb \lab{TanakaI}
    \Lambda^{x, \epsilon}_{t} := \int_0^t \<(-G_1 + \lambda)Q^{\lambda}_{ \epsilon}(x - \cdot), \mu_s\>ds. 
    \eeqlb 
    
  \bth \label{ClaimAk}
 Assume Hypotheses $\ref{hyp:basicassumpFilter}$, $\ref{hyp:basicassumpElliptic}$, 
 $\ref{hyp:basicassumpGauss}$ and $\ref{hyp:basicassumpUniformInteg}$.
  For any real numbers $T>0$, $\epsilon>0$ and any integer $k\ge1$, there holds
   \beqlb \lab{AME}
 \sup_{x \in \R^d} \sup_{0 \leq t \leq T}\E_{\mu_0} [|\Lambda^{x, \epsilon}_t |^k] < \infty
 \eeqlb
and
\beqlb \lab{BME}
   \lim_{\epsilon, \epsilon^{'} \downarrow 0} \sup_{x \in \R^d} \sup_{0 \leq t \leq T} \E_{\mu_0}
   [|\Lambda^{x, \epsilon}_t - \Lambda^{x, \epsilon{'}}_t |^k] = 0.
\eeqlb
 \eth
   
Since each $\Lambda^{x, \epsilon}_t(\omega) : [0, \infty) \times \R^d \times \Omega \rightarrow[0, \infty)$ is a jointly measurable function, 
there exists a common $\P_{\mu_0}$-null set $N$ for any $(t,x) \in [0, \infty) \times \R^d$, such that
\beqlb \lab{DELT}
 \Lambda^x_t := \begin{cases}
     \lim_{\epsilon \downarrow 0}\Lambda^{x, \epsilon}_t , & \quad \mbox{ if $\omega \notin N $} \\
     0 ,                                                  & \quad \mbox{ if $\omega \in N $}
     \end{cases}
\eeqlb
is a well-defined limit (in the $L^k(\P_{\mu_0})$ sense for every integer $k\ge1$) satisfying both 
\beqlb \lab{AME1}
 \sup_{x \in \R^d} \sup_{0 \leq t \leq T} \E_{\mu_{0}} [|\Lambda^{x}_t|^k] < \infty \nonumber   
\eeqlb
and
\beqlb \lab{BME1}
 \lim_{\epsilon \downarrow 0} \sup_{x \in \R^d} \sup_{0 \leq t \leq T} \E_{\mu_0}[|\Lambda^x_t - \Lambda^{x, \epsilon}_t|^k] = 0. \nonumber 
\eeqlb

The proof of Theorem \ref{ClaimAk} is achieved in Section \ref{sec:dualConst} by computing the $k^{th}$-moments of SDSM through a duality argument.

We now state our main result, which confirms that $(\ref{DELT})$ defines the local time of SDSM and possesses a Tanaka representation, 
under some restriction on the family of initial measures.


  \bth \lab{lt_th1}
 Under Hypotheses $\ref{hyp:basicassumpFilter}$ and $\ref{hyp:basicassumpElliptic}$, with $d=1$, $2$ or $3$,
 select any $a\ge0$ and ${\mu}_{0}\in M_a(\R^d)$ satisfying both 
 Hypotheses $\ref{hyp:basicassumpGauss}$ and $\ref{hyp:basicassumpUniformInteg}$, 
 with the Brownian sheet $W$ provided by Hypothesis $\ref{hyp:basicassumpFilter}$
 and martingale measure $M$ constructed from Theorem $\ref{MPforSDSM}$. Fix $\lambda>0$. 
 Outside of a $\P_{\mu_0}$-null set which is common to all $(t,x) \in (0, \infty) \times \R^d$, 
  $\Lambda^x_{t}$ defined by $(\ref{DELT})$ satisfies $\Lambda^x_{t} \geq 0$ and the following Tanaka formula holds: 
   \beqlb \lab{TanakaII}
    \Lambda^{x}_{t} & = & \<Q^{\lambda}(x - \cdot), \mu_0\> - \<Q^{\lambda}(x- \cdot), \mu_{t}\>
    + \lambda \int_0^{t} \<Q^{\lambda}(x - \cdot), \mu_s\>ds \nonumber \\
    & &  + \sum_{p=1}^{d}\int_0^{t} \int_{\R^d}\<h_p(y- \cdot)
    \partial_p Q^{\lambda}(x - \cdot), \mu_s\>W(dy,ds) \nonumber \\
   & &  + \int_0^{t} \int_{\R^d}Q^{\lambda}(x-y)M(dy,ds).
   \eeqlb
  Further, $\Lambda^x_{t}$ is the local time for SDSM $\{\mu_t\}$, in that 
  $(\ref{LT})$ holds outside of a $\P_{\mu_0}$-null set which is common to all $(t,\phi) \in [0, \infty) \times C_c(\R^d)$. 
 Finally, for every $t \ge 0$ there holds
 $ \sup_{x\in\R^d}\sup_{0 \leq s \leq t} \E_{\mu_0} [( \Lambda^x_{s} )^k] < \infty$ for every integer $k\ge1$.
 \eth
  \noindent  {\bf Remark}
Note that the $\P_{\mu_0}$-null set is independent of the test function $\phi \in {\cal S}(\R^d)$ in Theorem \ref{MPforSDSM} 
and of $x\in\R^d$ in Theorem \ref{lt_th1}, hence the same holds also in $(\ref{LT})$. 
Note also that the value of the local time does not depend on parameter $\lambda > 0$ 
(although it does vary with the dimension $d$ of the space). 

Theorem \ref{lt_th1} establishes the existence of the local time $\Lambda^x_t$ for SDSM directly through the characterization provided by (\ref{TanakaI}) and (\ref{DELT}), an explicit Tanaka formula (\ref{TanakaII}) expressed through a Green function with a spatial singularity, 
in the spirit of the approach proposed in L\'opez-Mimbela and Villa \cite{Lopez-MimbelaVilla04} in their Theorem 3.1, 
of which our Theorem \ref{lt_th1} is an extension.
However, in order to make sense of it, we have to approximate this singular Green function
and its derivatives by smooth functions to ensure that the various stochastic integrals in (\ref{TanakaII}) are well-defined. 

Recently, using Malliavin Calculus, Hu et al. \cite{HNX19} proved H{\"o}lder continuity for a
related class of processes, namely SDSM with the classical Dawson-Watanabe branching
mechanism replaced by the more cohesive and therefore asymptotically smoother Mytnik-Sturm branching,
under an initial measure $\tilde\mu_0$ with a bounded Radon-Nikodym derivative with respect
to Lebesgue measure $\lambda_0$. More precisely they have shown that, in the case where
$c$ is the identity matrix and $h$ is smooth enough (and matrix-valued),
this regularized SDSM $\{\tilde\mu_t\}$ has a density $\tilde f_t=d\tilde\mu_t/d\lambda_0$
which is almost surely jointly H{\"o}lder continuous,
with exponent $\beta_1\in(0,1)$ in space and $\beta_2\in(0,1/2)$ in time. In fact they showed that for every such $\beta_1$
and $\beta_2$, as well as for every choice of $p > 1$, there is a constant $c=c(T, d, h, p, \beta_1, \beta_2)> 0$ such that,
for every $x,z\in\R^d$ and $0 < s < t \le T$ there holds
\[
 \left[ \E_{\tilde\mu_0} \left| \tilde f_t(z)-\tilde f_s(x)\right|^{2p} \right]^{1/2p}
 \leq cs^{-1/2} (|z - x|^{\beta_1}+|t - s|^{\beta_2}).
\]
Under these stronger initial conditions we can write the regularized local time as
\[
 \tilde\Lambda^x_{t} = \int_0^t \tilde f_s(x)ds
 \]
which is thus almost surely jointly H{\"o}lder continuous, with (at least) the same exponents.
Such a phenomenon is not likely to occur for the local time of our SDSM here, since there is such a density
 $f_t=d\mu_t/d\lambda_0$ when $d=1$ (Konno and Shiga \cite{KonnoShiga88} for Super-Brownian motion
 and Dawson et al \cite{DVW2000} for the general case)
 but not when $d\ge2$, not even in the Super-Brownian motion case
 (see Dawson and Hochberg \cite{DawsonHochberg79} and Perkins \cite{Perkins88}, \cite{Perkins02}).
 The sharpest estimates for the closed support of $\mu_t$  in the Super-Brownian motion case are found in
 Dawson and Perkins \cite{DawsonPerkins99} when $d\ge3$ and
 Le Gall and Perkins \cite{LeGallPerkins95} when $d=2$. See also Hong \cite{Hong18} for renormalization issues.
 The matter of the modulus of continuity of the local time $\Lambda^x_{t}$,
 including in the aforementioned special cases, is addressed in a separate paper \cite{DVW2021}.


\section{Branching model, SDSM and SPDE}\label{sec:SDSM}
\setcounter{equation}{0} 

The construction of SDSM as a limit of a sequence of branching particle systems in $\R^d$ is presented here. 
First we consider the simplest situation in which there is a finite number $m$ of particles moving in $\R^d$ without branching, 
but submitted to a diffusive dynamic that includes an interacting property, which is produced by the random environment. 
For this we need some useful properties of finite dimensional diffusions.

\subsection{Diffusive part of the finite systems and Hypothesis \ref{hyp:basicassumpGauss}}

Consider a system of $m\ge1$ interacting particles $\{z_{k}: k=1,\ldots,m\}$, each moving continuously in
$\R^d$  from some arbitrary starting point $z_{k}(0)\in \R^d$ and each driven by its own standard $d$-dimensional
Brownian motion $B_{k}$ independently of one another, but all of them evolving coherently in a random medium prescribed by a
common Brownian sheet $W$ on $\R^d$, which remains the same throughout the construction and is assumed to be independent of
the set $\{B_{k}: k\ge1\}$. Under Hypotheses $\ref{hyp:basicassumpFilter}$ and $\ref{hyp:basicassumpElliptic}$   
on diffusion coefficient $c$ (a $d\times d$ matrix of real-valued functions)
and random medium intensity $h$ (a vector of $d$ real-valued functions), the system of stochastic integral equations 
  \beqlb \label{individual} z_k(t)=z_k(0) + \int_{0}^t c(z_{k}(s)) dB_{k}(s)
 + \int_{0}^{t } \int_{\R^d}h(y-z_{k}(s))W(dy,ds)
 \eeqlb
for $k=1,\ldots,m$, has a unique strong solution for any $m\ge1$ 
and every fixed starting point $(z_1(0), z_2(0), \ldots, z_m(0))\in (\R^d)^m$. 
It is a strong $\{{\cal F}_t\}_{t \geq 0}$-Markov process and there is no explosion in finite time.
Existence and uniqueness of solution follow by Picard's method of successive approximations,
 as in Wang \cite{Wang97}. Bounded continuous coefficients preclude explosion. 
In other words (\ref{individual}) has a strong solution and pathwise uniqueness holds, in the sense that
any two solutions with sample paths $\P$-almost surely in $C([0, \infty), (\R^d)^m)$ and identical starting points,
must be equal with $\P$-probability $1$.
The cloud of particles $\{z_{1}, z_{2}, \ldots , z_{m}\}$ as a whole clearly diffuses according to correlated dynamics, an important 
feature of the motion is that it generates new difficulties in the identification of the local time of SDSM, as will be seen shortly.

For any integer $m \geq 1$, write $Z_m(t):= (z_{1}(t), \cdots, z_{m}(t))$
for the motion of the cloud of $m$-particles solving (\ref{individual}), $\P_x$ for the law of $Z_m$ with initial
point $x\in (\R^d)^m$ and $\E_x$ for the expectation with respect to $\P_x$.
Since $Z_m$ is a time-homogeneous $\{{\cal F}_t\}_{t \geq 0}$-Markov process, let $\{P^m_t: t\ge 0\}$ be the
corresponding Markov semigroup on $B({(\R^d)^m})$ for $Z_m$, that is
 \beqlb\label{eqn:Semigroup}
 P^{m}_t f(x):=\E_x \left[ f(Z_{m}(t) )\right] \qquad
\hbox{for } t\geq 0 \hbox{ and } f \in B({(\R^d)^m}).
 \eeqlb
Note that $P^m_t$ is a Feller semigroup and maps each of $B({(\R^d)^m})$,
$C_b({(\R^d)^m})$ and $C_0({(\R^d)^m})$ into itself.

 It\^{o}'s formula yields the following generator for $\{P^m_t: t\ge 0\}$ : for all $f\in C_b^{2}({(\R^d)^m})$, 
 \beqlb\label{eqn:Gn}
{G}_{m} f(x):=
   \frac{1}{2} \sum_{i,j=1}^{m} \hspace{1mm}
 \sum_{p, q=1}^{d}
 \Gamma_{pq}^{ij}(x_{1}, \cdots, x_{m})
  \frac{\partial^2}{\partial
 x_{ip} \partial
 x_{jq} }
 f(x_{1}, \cdots, x_{m})
 \eeqlb
 where $x=(x_1, \cdots, x_m)\in (\R^d)^m$ has components
 $x_i = (x_{i1}, \cdots, x_{id}) \in \R^d$ for $1\leq i \leq m$ and
 $\Gamma_{pq}^{ij}$ is defined by $(\ref{gammaij})$.
 Hypothesis $\ref{hyp:basicassumpElliptic}$ ensures that operator ${G}_{m}$ is uniformly elliptic --- 
 see Subsection \ref{app:pf_ellip1} for a proof. 
 
Following Stroock and Varadhan \cite{StroockVaradhan79}, it is useful to view process $\{Z_m(t):{t \geq 0}\}$
as a solution to the $({G}_{m}, \delta_{Z_m(0)})$-martingale problem on 
$(\Omega, {\cal F}, \{{\cal F}_t\}_{t \geq 0}, \P)$ for any fixed starting point $Z_m(0)\in (\R^d)^m$, meaning that, for every 
choice of $f\in C_c^{\infty}({(\R^d)^m})$, the process $f(Z_m(t))-\int_{0}^{t} {G}_{m} f(Z_m(s))ds$ is an ${\cal F}_t$-martingale.
We say this martingale problem is well-posed (or has a unique solution) if there exists at least one solution and any two solutions 
have the same finite dimensional distributions. 

We also need the following summary of several known results from the literature.

\ble \lab{FellerProp}
Under Hypotheses $\ref{hyp:basicassumpFilter}$ and $\ref{hyp:basicassumpElliptic}$,
the following statements hold for every choice of $m\geq 1$.
\begin{itemize}
\item
For any initial value $Z_m(0)\in (\R^d)^m$, the $({G}_{m}, \delta_{Z_m(0)})$-martingale problem is well-posed.
The trajectories of $\{Z_m(t):{t \geq 0}\}$ are in $C([0, \infty), (\R^d)^m)$.
\item
$P_t^m f (x)$, as a function of $(t,x)$, belongs to $\cup_{t>0}C_b^{1,2}((0, t] \times (\R^d)^m)$ when $t>0$,
for every choice of $f\in C_0({(\R^d)^m})$, and $\{P_t^m\}$ is a Feller semigroup mapping
$C_0^2((\R^d)^m)$ into itself.
\item
$\{P_t^m: t\ge 0\}$ has a transition probability density when $t > 0$, i.e., there is
a function $q_t^{m}(x,y) > 0$ which is jointly continuous
in $(t,x,y) \in (0, \infty) \times (\R^d)^m \times (\R^d)^m$
everywhere and such that there holds
$P_t^m f (\cdot)=\int_{(\R^d)^m} f(y) q_t^{m}(\cdot,y)dy$
when $t > 0$, for every $f\in C_0({(\R^d)^m})$.
\item
For each choice of $T > 0$, $d\ge1$ and $m\ge1$, there are positive constants $a_1$ and $a_2$ such that,
for any choice of $1\le p,p' \le dm$ and nonnegative integers $r$, $s$ and $s'$ verifying $0\le 2r+s+s'\le2$,
\beqlb \lab{LSU}
   \left|\frac{\partial^{r}}{\partial t}\frac{\partial^{s}}{\partial y_{p}}\frac{\partial^{s'}}{\partial y_{p'}}  q_t^{m}(x,y)\right|
   \leq \frac{a_1}{t^{(dm+2r+s+s')/2}}\exp{\left\{- a_2\left( \frac{|y-x|^{2}}{t} \right)\right\}}
  \eeqlb
 holds everywhere in $(t,x,y)\in(0,T) \times (\R^d)^m \times (\R^d)^m$ with $y=(y_1,\ldots,y_{dm})$.
  \item
For each choice of $T > 0$, $d\ge1$ and initial data $(\xi, \tau) \in \R^d \times [0,T)$, 
there is a unique fundamental solution $\Gamma(x,t; \xi, \tau)$ to ${G}_{1}\Gamma-\partial_t\Gamma = 0$.
Moreover, there exist two constants $c > 0$ and $c_0 > 0$, such that,
for all nonnegative integers $r,s_1,s_2,\ldots,s_d$ verifying
$0\le l=2r + s_1+s_2+ \ldots + s_d \le 2$ and writing
$\partial^{l}=\partial^{r}_t\partial^{s_1}_{x_1}\partial^{s_2}_{x_2}\cdots\partial^{s_d}_{x_d}$
with $\partial^{0}$ for the identity, there holds, for every choice of $\alpha\in(0,1)$,
 \beqlb \lab{formulaI}
 & & |\partial^{l}\Gamma(x, \xi; t, \tau) -\partial^{l}\Gamma(x, \xi^{'}; t, \tau^{'})|  \nonumber  \\
\leq & & c( |\xi-\xi^{'}|^{\alpha} + |\tau - \tau^{'}|^{\alpha/2})
\bigg[ (t-\tau)^{- (d+l)/2} \exp{ \{- c_0 \frac{|x - \xi|^2}{t - \tau} \} } \nonumber  \\
 & & \hspace{5cm} + (t-\tau^{'})^{- (d+l)/2} \exp{\{-c_0 \frac{|x - \xi^{'}|^2}{t - \tau^{'}}\}}  \bigg].
 \eeqlb
  \end{itemize}
 \ele

 \noindent  {\bf Remark}
 One important consequence of Lemma \ref{FellerProp} is that
 $C_0^2((\R^d)^m)$ is a core for generator ${G}_{m}$
 (see Propositions 1.3.3 and 8.1.6 in Ethier and Kurtz \cite{EthierKurtz86}).

\proof
Since $h_p \in L^1(\R^d) \cap {\rm Lip}_b (\R^d) \subset L^2(\r^d)$ holds for all $p\ge1$,
so does $\rho_{pq}\in {\rm Lip}_b(\R^d)$
and $\Gamma_{p,q}^{i,j}\in {\rm Lip}_b((\R^d)^m)$ for all $p,q\ge1$.
Hypothesis \ref{hyp:basicassumpElliptic} is stronger than the conditions in each of Theorem 6.3.4 p.152,
Theorem 3.2.1 p.71 or Corollary 3.2.2 p.72 in Stroock and Varadhan \cite{StroockVaradhan79}, hence
the Feller property holds true for Markov semigroup $P^m_t$ in (\ref{eqn:Semigroup})
and the first three statements follow, using the uniform density of $C_c^{\infty}({(\R^d)^m})$ in $C_0({(\R^d)^m})$.
The fourth, the upper bound in (\ref{LSU}) on the transition density $q_t^{m}(x,y)$ of semigroup $\{P_t^{m}: t\ge 0\}$
generated by $G_{m}$ from (\ref{eqn:Gn}), is a consequence of Equation (13.1) of Lady\u{z}enskaja et al. (\cite{LSU68} p.376). 
The fifth ensues from Theorem 5.3.5 of Garroni and Menaldi \cite{GarroniMenaldi92}, 
as ${G}_{1}-\partial_t$ is uniformly parabolic whenever ${G}_{1}$ is uniformly elliptic. 
\qed

 \noindent  {\bf Remark:} 
When combined with Hypothesis $\ref{hyp:basicassumpElliptic}$, condition $(\ref{nonBA})$ 
in Hypothesis \ref{hyp:basicassumpGauss} is equivalent to 
 \beqlb \lab{H4}
 \sup_{x\in \R^d}  \sup_{0 < t \leq T} \<q_t^{1}(\cdot,x), \mu_0\>  < \infty,
 \eeqlb
where $q_t^{1}$ is the one particle transition density with generator $G_1$ on $\R^d$ 
given by either $(\ref{eqn:Gn1})$ or $(\ref{eqn:Gn})$. This is a consequence of the uniform ellipticity,
which provides us with both an upper bound $(\ref{LSU})$ as well as the
corresponding lower bound : there exist four positive constants $a^*$, $b$, $c$ and $A^*$ such that
\bea \label{Aronsonbounds} 
   a^* \cdot  \varphi_{bs} (y-x) \leq q_s^{1}(x,y) \leq A^* \cdot \varphi_{cs}(y-x)
\eea
holds for any $x, y \in \R^d$ and $s > 0$ (Aronson \cite[Theorem 10]{Aronson68}). 
Notice also that any measure $\mu_0$ with a finite Radon-Nikodym derivative
with respect to Lebesgue measure $\lambda_0$ (finite as in finitely $\lambda_0$-integrable) satisfies
$ \sup_{t > 0} \sup_{x\in \R^d} \< \varphi_t(x-\cdot), \mu_0\>  < \infty$ and therefore
Hypothesis $\ref{hyp:basicassumpGauss}$ as well. In particular this is the case for measures
$I_a(x)dx$, for all choices of $a\ge0$. 

The need for Hypothesis $\ref{hyp:basicassumpGauss}$ stems from the following illustration. 
 
 \begin{example} Let $q_t^{1}$ is the one particle transition density with generator $G_1$ on $\R^d$,
 made explicit in (\ref{eqn:Gn}) above.
For initial measure $\mu_0=\delta_0$ and any $t >0$, there holds
 \bea \label{ex} \hspace*{8mm}\int_0^t \int_{\r^d}q_s^{1}(y,x) \delta_0(dy) ds
 = \int_0^t q_s^{1}(0,x)ds  \left\{
\begin{array}{lll}
                                < \infty  \quad    & \mbox{if  $x=0$},  &  d=1 \\
                                 < \infty \quad   &  \mbox{if  $x \neq 0$}        &  d=1, \\
                                 = \infty \quad    & \mbox{if  $x_1x_2=0$},  &  d=2 \\
                                 < \infty \quad   &  \mbox{if  $x_1x_2 \neq 0$}        &  d=2.
                                 \end{array}
                                 \right.
\eea
\end{example}
From this example, we see that, if the initial measure has an atom and $d \geq 2$,
the existence of a continuous local time for SDSM is questionable.
This motivates the constraint on the family of initial measures set forth in Hypothesis \ref{hyp:basicassumpGauss}. 

 \noindent  {\bf Remark:} 
 Lebesgue measure $\lambda_0$ on $\R^d$ satisfies Hypothesis \ref{hyp:basicassumpGauss} for all $d\ge1$ but 
 satisfies Hypothesis \ref{hyp:basicassumpUniformInteg} when and only when $a>d$. In general, when $a>d$, 
 any measure $\mu_0 \in M_a(\R^d)$ which satisfies Hypothesis \ref{hyp:basicassumpGauss} will also satisfy  
 Hypothesis \ref{hyp:basicassumpUniformInteg}. Indeed, using  
 $\<I_a(\cdot+w), \mu_0 \> = \lim_{t\downarrow0}\int_{\R^d}\int_{\R^d}I_a(x+w+y)\varphi_{t}(y)\lambda_0(dy)\mu_0(dx)$, 
 the translation invariance of Lebesgue measure yields  
 
  \beqlb \lab{Tonelli} 
   \int_{\R^d}\int_{\R^d}I_a(x+w+y)\varphi_{t}(y)\lambda_0(dy)\mu_0(dx) 
  &=&\int_{\R^d}\int_{\R^d}I_a(w+y)\varphi_{t}(y-x)\lambda_0(dy)\mu_0(dx) \nonumber \\
  &\le& \sup_{y\in \R^d} \sup_{0 < t \leq T} \<\varphi_{t}(y-\cdot), \mu_0\> 
\cdot \sup_{w\in \R^d} \<I_a(w+\cdot), \lambda_0\>  \nonumber \\ 
  && < \infty  
\eeqlb 
and $\sup_{w\in \R^d} \<I_a(\cdot+w), \mu_0 \> < \infty$ ensues.


\subsection{Branching Particle Systems}

Ren et al. \cite{RenSongWang09} handled the construction of SDSM and the derivation of the associated SPDE,
 using a tightness argument for the laws on $D([0, \infty), M_0(D))$ of the trajectories of high-density particles,
 but only when these particles move in a bounded domain $D \subset \R^d$
 with killing boundary and the initial data is a finite measure $\mu_{0} \in M_0(D)$. 
There are some significant differences between the construction of our model in $\R^d$ and theirs, 
notably in the topological structure of the spaces required, 
so, for the sake of completeness, some details are provided here about the 
relevant branching particle systems in $\R^d$ when $d\ge1$, 
as well as the construction of SDSM in $\R^d$ and the associated SPDE on $\R^d$. 
Additional details can be found in Ren et al. \cite{RenSongWang09}. 

A prior construction for SDSM on the whole real line ($d=1$) can be found in Dawson et al. \cite{DLW01}, where 
the initial data is bounded and use is made of the canonical space $C([0, \infty), M_0(\R))$ since the trajectories are 
less irregular than in the current case ($d\ge2$), thus affording sharper inequalities. 
Earlier still, Konno and Shiga \cite{KonnoShiga88} built super stable processes on $\R^d$ with unbounded initial data 
but no interacting term ($h\equiv0$). In both cases the construction is on the canonical space $C([0, \infty), M_a(\R^d))$, 
an approach that is unavailable here since the interacting term depends on space-time stochastic integrals that control 
the evolution of the strong solution of the associated SDE. 

Therefore we must also begin with the full definition of the interacting branching particle systems. 

In summary, each particle has an exponentially distributed lifetime at the end of which it either splits into finitely many identical particles or
dies, independently of one another and of the above diffusion mechanism. The newborn particles (if there are any) then start afresh, 
with their own independent exponential lifetimes, at the spatial position where the parent left off and the new (enlarged or diminished) 
cloud of particles continues its collective evolution in $\R^d$ according to (\ref{individual}). 
For the resulting branching particle systems,
particles undergo a finite-variance branching at independent
exponential times and have interacting spatial motions powered by
diffusions and a common white noise. 
Under critical conditions on the branching rate, the distribution of the number of offsprings 
(the mean number of offspring is asymptotically one) and the (common) mass of the particles 
ensuring a non trivial high density limit, the sequence
of empirical measure processes representing the proportion of the mass carried by those particles alive at a given time and located
within a given set, converges to a limiting measure-valued Markov process $\{\mu_t:t\ge0\}$. 

 We first introduce an index set in order to identify each particle in the branching
 tree structure. Let $\Re$ be the set of all multi-indices, i.e.,
 strings of the form
$\xi = n_{1} \oplus n_{2} \oplus \cdots \oplus n_{k}$, where the
$n_{i}$'s are non-negative integers. Let $| \xi |$ denote the length
of $\xi$. We provide $\Re$ with the arboreal ordering: $m_{1} \oplus
m_{2} \oplus \cdots \oplus m_{p} \prec n_{1} \oplus n_{2} \oplus
\cdots \oplus n_{q}$ if and only if $p \leq q$ and $m_{1}=n_{1},
\cdots, m_{p}=n_{p}.$ If $| \xi | = p$, then $\xi$ has exactly $p-1$
predecessors, which we shall denote respectively by $\xi - 1$, $\xi
- 2, \cdots, \xi - | \xi |+1 $. For example, with $\xi = 6 \oplus 18
\oplus 7 \oplus 9$, we get $ \xi - 1 = 6 \oplus 18 \oplus 7$, $\xi -
2 = 6 \oplus 18$ and $\xi - 3 = 6$. We also define an $\oplus$
operation on $\Re$ as follows: if $\eta \in \Re$ and $|\eta| = m$,
for any given non-negative integer $k$, $\eta \oplus k \in \Re$ and
$\eta \oplus k$ is an index for a particle in the $(m + 1)$-th
generation. For example, when $\eta = 3 \oplus 8 \oplus 17 \oplus 2$
and $k = 1$, we have $\eta \oplus k = 3 \oplus 8 \oplus 17 \oplus 2
\oplus 1 $.

Let $\{B_{\xi}=(B_{\xi 1}, \cdots, B_{\xi d})^T: \, \xi \in \Re\}$
be an independent family of standard $\R^d$-valued Brownian motions,
where $B_{\xi k}$ is the $k$-th component of the $d$-dimensional
Brownian motion $B_{\xi}$, and $W$ a Brownian sheet on $\R^d$.
Assume that $W$ and $\{B_{\xi}: \xi \in \Re\}$ are defined on a
common filtered probability space $(\Omega, {\cal F}, \{{\cal
F}_t\}_{t\geq 0}, \P)$, and independent of each other. 
Under Hypotheses $\ref{hyp:basicassumpFilter}$ and $\ref{hyp:basicassumpElliptic}$, 
just as for equations (\ref{individual}) (see Lemma 3.1 of Dawson et al. \cite{DLW01}), for every
index $\xi\in \Re$ and initial data $z_{\xi}(0)$, by Picard's iteration method, there is a unique strong
solution $z_{\xi}(t)$ to the equation
 \beqlb \label{3.1}
 z^{T}_{\xi}(t)=z^{T}_{\xi}(0) +
\int_{0}^{t} c(z_{\xi}(s)) dB_{\xi}(s) +
\int_{0}^{t}\int_{\R^d}h(y-z_{\xi}(s))W(dy,ds).
 \eeqlb
Since the strong solution of (\ref{3.1}) only depends on the initial
state $z_{\xi}(0)$, the Brownian motion
$B_{\xi}:=\{B_{\xi}(t):t\ge0\}$ and the common $W$, this solution can be written as 
$z_{\xi}(t)=\Phi (z_{\xi}(0),B_{\xi},t)$ for some measurable $\R^d$-valued map $\Phi$  
(dropping $W$ for the sake of simplifying the formulas coming up). 
For every $\phi\in C_b^{2}(\R^d)$ and $t>0$, It\^{o}'s formula yields
 \beqlb \lab{3.4}
 && \hspace{-1cm}  \phi(z_{\xi}(t)) - \phi(z_{\xi}(0))   \\
  &=&  \sum_{p=1}^{d}\left[\int_{0}^{t}
 \left( \partial_{p} \phi(z_{\xi}(s)) \right)
 \sum_{i=1}^{d} c_{pi}(z_{\xi}(s))
 dB_{\xi i}(s) \right. \nonumber \\
  & & \left.+ \int_{0}^{t}
 \int_{\r^d} \partial_{p} {\phi}(z_{\xi}(s))
 h_p(y- z_{\xi}(s))W(dy,ds) \right] \nonumber \\
  & & +
 \frac{1}{2} \sum_{p,q=1}^{d}\int^{t}_{0}\left( \partial_{p}
 \partial_{q} {\phi}(z_{\xi}(s))
 \sum_{i=1}^{d}c_{pi}(z_{\xi}(s)) c_{qi}(z_{\xi}(s))
\right) \,
 ds    \nonumber \\
& & + \frac{1}{2} \sum_{p,q=1}^{d} \int_{0}^{t}(
\partial_{p}
\partial_{q} {\phi}(z_{\xi}(s))) \int_{\R^d} h_p(y- z_{\xi}(s) )
 h_q(y- z_{\xi}(s) )dy ds.
    \nonumber
 \eeqlb

With the spatial motion of each particle modelled by (\ref{3.1}), 
we next spell out the critical conditions governing the branching mechanism. 
For every positive integer $n\geq 1$, there is an initial system of
$m_{0}^{(n)}$ particles. Each particle has mass $1/ {\theta}^{n}$
and branches independently at rate $\ga \theta^{n}$. Let $q^{(n)}_k$
denote the probability of having $k$ offspring when a particle survives
in $\R^d$. The sequence $\{q^{(n)}_{k}\}$ is assumed to satisfy the
following conditions:
$$
q^{(n)}_{k} = 0 \qquad \hbox{if } k=1 \hbox{ or } k \geq n+1,
$$
and
$$
\sum_{k=0}^{n} k q^{(n)}_{k} =1 \quad \hbox{ and } \quad \lim_{n \ra
\infty}\sup_{k\geq 0} |q^{(n)}_{k}-p_{k}|=0 ,
$$
where $\{p_k: k=0, 1, 2, \cdots \}$ is the limiting offspring
distribution which is assumed to satisfy following conditions:
$$
p_{1}=0 , \hspace{1mm} \sum_{k=0}^{\infty} kp_{k} =1 \hspace{1mm}
\hbox{ and } \hspace{1mm}  m_2 := \sum_{k=0}^{\infty}k^{2}p_{k} <
\infty.
$$
Let $m_c^{(n)}:= \sum_{k=0}^n(k-1)^4q_k^{(n)}$. The sequence
$\{m_c^{(n)}: n\geq 1\}$ may be unbounded,  but we assume that
$$
\lim_{n \ra \infty}\frac{m_c^{(n)}}{\theta^{2n}}=0
\hspace{1cm}\mbox{ for any $\theta > 1$. }
$$
Assume that $m_{0}^{(n)} \leq \hbar \, {\theta}^{n} $, where $\hbar> 0$ and $\theta > 1$ are fixed constants. 
Finally, define
$m^{(n)}_2 := \sum_{k=0}^{n} k^{2}q^{(n)}_{k}$, 
$\sigma^2_{n} :=m^{(n)}_2 -1$ and $\sigma^2 := m_2 -1$. Note that $\sigma^2_n$ and
$\sigma^2$ are the variance of the $n$-th stage and the limiting
offspring distribution, respectively. We have $\sigma^2_{n} <
\infty$ and $\lim_{n \ra \infty} \sigma^2_n = \sigma^2$.

 For a fixed stage $n \ge 1$, particle $\xi \in\Re$ is located at $x_{\xi}(t)$ at time $t$; 
 $\{O^{(n)}_{\xi}:\,  \xi\in\Re\}$ is a family of i.i.d. random variables with
$\P(O^{(n)}_{\xi} = k) = q^{(n)}_{k}$ and $k=0,1,2,\cdots$; and
$\{C^{(n)}_{\xi}: \, \xi \in \Re\}$ be a family of i.i.d.
real-valued exponential random variables with parameter $\ga
{\theta}^{n}$, which will serve as lifetimes of the particles. We
assume $W$, $\{B_{\xi}: \, \xi \in{\Re}\}$, $\{C^{(n)}_{\xi}: \, \xi
\in{\Re}\}$ and $\{O^{(n)}_{\xi}: \, \xi \in\Re\}$ are all
independent. Once the particle $\xi$ dies, it is sent at once to a cemetery point noted by $\partial$.

To simplify our notation further, we drop the superscript $(n)$ from the random variables. 
 If the death location of the particle $\xi-1$ is in $\R^d$,
then the birth time $\beta(\xi)$ of the particle $\xi$ is given by
$$
\beta(\xi) := \left\{ \ba{ll} \sum_{j=1}^{| \xi |-1}C_{\xi - j}, &
\mbox{if $ O_{\xi -j} \geq 2$ for every $j = 1, \cdots , |\xi|-1$\,;} \\
\infty,   & \mbox{otherwise}. \ea \right.
$$

The death time of the particle $\xi$ is given by $\zeta(\xi) =
\beta(\xi) + C_{\xi}$ and the indicator function of the lifespan of
$\xi$ is denoted by $\ell_{\xi}(t) := 1_{[\beta(\xi), \zeta(\xi)
)}(t)$.

Define $x_{\xi}(t) = \partial$ if either $ t < \beta(\xi) $ or $ t \geq
\zeta(\xi)$. We make the convention that any function $f$ defined on
$\R^d$ is automatically extended to $\R^d \cup \{
\partial \}$ by setting $f(\partial) = 0$.

Given $\mu_{0}\in { M}_{a}( \R^d)$, let $\mu_0^{(n)}:=(1/\theta^n)\sum_{\xi=1}^{m_{0}^{(n)}}\delta_{x_{\xi}(0)}$ 
be constructed upon a collection of initial starting points $\{x_{\xi}(0) \}$ for each $n\geq1$, in order for 
$\mu_0^{(n)} \Ra\mu_{0}$ to hold as $n \ra \infty$.

 Let ${\cal N}_{1}^{n} := \{1,2, \cdots, m_{0}^{(n)} \}$ be the set of indices
for the first generation of particles. For any $\xi \in{\cal
N}_{1}^{n} \cap \Re$, if $x_{\xi}(0) \in \R^d$, define \bea \label{3.5}
 x_{\xi}(t):= \begin{cases} \Phi (x_{\xi}(0), B_{\xi}, t),
&\qquad t \in [0,  C_{\xi}) ,  \\
\partial , &\qquad t < 0, \mbox{or }  t \geq C_{\xi},
            \end{cases}
\eea
and
$$
x_{\xi}(t) \equiv \partial \hspace{8mm}\mbox{for any $\xi \in
(\N \setminus {\cal N}_{1}^{n}) \cap \Re$ and $t \geq 0$}.
$$
For the path of the second generation, let $\bar{\zeta_1} = min\{C_{\xi}: \xi \in {\cal N}_{1}^{n} \cap \Re\}$. 

By Ikeda et al. \cite{INW68} and \cite{INW69}, for each $\omega \in \Omega$, there exists a measurable selection $\xi_0= \xi_0(\omega) \in {\cal
N}_{1}^{n} \cap \Re $ such that $\bar{\zeta_1} = C_{\xi_0}$. 
If $\xi_0\in {\cal N}_1^n \cap \Re$ and $x_{\xi_{0}}(t_0) = \partial$ for some $t_0 > 0$, then
$x_{\xi}(t) \equiv \partial$ for any $\xi \succ \xi_0 $ and any $t
\geq \zeta({\xi}_{0})\vee t_0$. Otherwise, if $x_{{\xi}_{0}}(\zeta({\xi}_{0})-) \in \R^d$ and
$O_{{\xi}_{0}}(\omega)=k \geq 2$, define for every ${\xi} \in
\{{\xi}_{0} \oplus i: \ i=1,2, \cdots,k \}$,
 \bea \label{3.6}
 x_{\xi}(t):=\begin{cases}
\Phi (x_{{\xi}_{0}}(\zeta({\xi}_{0})-), B_{\xi}, t), & \quad  t \in
[\beta(\xi) , \,  \zeta({\xi})), \\
\partial ,       &\quad  t \ge \zeta({\xi}).
\end{cases}
 \eea
 If $O_{{\xi}_{0}}(\omega)=0$, define $x_{\xi}(t) \equiv
\partial$ for $0 \leq t < \infty$ and ${\xi} \in \{{\xi}_{0} \oplus i: i
\geq 1 \}$.

More generally for any integer $m\geq 1$, let ${\cal N}_{m}^{n}
\subset \Re$ be the set of all indices for the particles in
 the $m$-th generation. If $\xi_0\in {\cal N}_m^n$ and if
$x_{{\xi}_{0}}(\zeta({\xi}_{0})-) \in \partial$, then $x_{\xi}(t)
\equiv \partial$ for any $\xi \succ \xi_0$ and any $t \geq
\zeta({\xi}_{0})
$. Otherwise, if $x_{{\xi}_{0}}(\zeta({\xi}_{0})-) \in \R^d$ and
$O_{\xi_{0}}(\omega)=k \geq 2$, define for ${\xi} \in \{{\xi}_{0}
\oplus i: i=1,2, \cdots,k \}$
 \bea \label{3.7}
 x_{\xi}(t):=\begin{cases}
\Phi (x_{{\xi}_{0}}(\zeta({\xi}_{0})-), B_{\xi}, t), & \quad  t \in
[\beta(\xi) , \,  \zeta({\xi})), \\
\partial ,       &\quad  t \ge \zeta({\xi}).
\end{cases}
 \eea
  If $O_{{\xi}_{0}}(\omega)=0$, define
$x_{\xi}(t) \equiv \partial$ for $0 \leq t < \infty$ 
and  for ${\xi} \in \{{\xi}_{0} \oplus i: i \geq 1 \}$. 
Continuing in this way, we obtain a branching tree of particles for
any given $\omega$ with random initial state taking values in
$\left\{x_{1}(0),x_{2}(0),\cdots,x_{m_{0}^{(n)}}(0)\right\}$. 


\subsection{Tightness and SPDE for SDSM}

We proceed to show that the sequence of laws of the empirical processes 
  \beq \label{4.1}
\mu^{ (n)}_{t} := \frac{1}{{\theta}^{n}}\sum_{\xi \in{\Re}}\delta_{x_{\xi}(t)},
 \eeq
 associated with the branching particle system $\{x_{\xi}\}$ constructed in the last subsection, 
 converges weakly as $n\to \infty$ on $D([0, \infty), {\cal S}^{\prime}(\R^d))$, that its weak limit 
 is actually on $C([0, \infty), {\cal S}^{\prime}(\R^d))$ and that this limit is our SDSM from (\ref{SPDEc}). 
 
First recall that every finite positive Radon measure on $\R^d$ defines a tempered distribution on $\R^d$; 
the injection of $M_0(\R^d)$ into ${\cal S}' (\R^d)$ is in fact continuous by Proposition 5.1 in Dawson and Vaillancourt \cite{DV95}. 
We can therefore extend the bracket notation $\<\phi,\mu\>$ to the duality between ${\cal S} (\R^d)$ and ${\cal S}' (\R^d)$ without 
ambiguity. Sometimes we write $\mu(\phi)$=$\<\phi,\mu\>$. Furthermore this induces continuous injections of 
$C([0, \infty), M_0(\R^d))$ into $C([0, \infty), {\cal S}^{\prime}(\R^d))$ and 
$D([0, \infty), M_0(\R^d))$ into $D([0, \infty), {\cal S}^{\prime}(\R^d))$. 

We now build the transfer function that maps the infinite measure-valued processes to finite ones. 
For each $a\ge0$, the map $T_{I_a} : M_a(\R^d) \rightarrow M_0(\R^d)$
defined by
\beqlb \label{transfer}
 T_{I_a}(\mu)(A) := \int_{A} I_a(x) \mu(dx)=\int_{A} (1 + |x|^2)^{- a/2} \mu(dx),
\eeqlb
for any $A \in {\cal B}(\R^d)$, is clearly homeomorphic (continuous and bijective,
with a continuous inverse), hence also Borel measurable. It also induces continuous mappings from 
$C([0, \infty), M_a(\R^d))$ into $C([0, \infty), {\cal S}^{\prime}(\R^d))$ and 
$D([0, \infty), M_a(\R^d))$ into $D([0, \infty), {\cal S}^{\prime}(\R^d))$. 

For any $t>0$ and $A \in {\cal B}({\R}^d)$,  define
\beq \label{4.2}
 M^{(n)}(A \times (0,t]) := \sum_{\xi \in {\Re}} \frac{[O^{(n)}_{\xi}
 -1]}{{\theta}^{n}}
 1_{\{x_{\xi}({\zeta(\xi)-})
 \in A, \, \zeta(\xi) \leq t\}}, \,
\eeq
 which describes the space-time related branching in the set $A$ up to time $t$. 

By (\ref{3.4}) we have a finite dimensional SPDE
 \beqlb \lab{DSPDE}
  & &\hspace{-1cm} \<\phi(\cdot) I_a^{-1}(\cdot), I_a(\cdot)\mu^{(n)}_{t}\> -\<\phi(\cdot) I_a^{-1}(\cdot), I_a(\cdot) \mu^{(n)}_{0}\> 
  \nonumber \\
  & = & \frac{1}{\sqrt{{\theta}^{n}}}\<\phi(\cdot) I_a^{-1}(\cdot), I_a(\cdot)U^{(n)}_{t}\> + \<\phi(\cdot) I_a^{-1}(\cdot), I_a(\cdot)X^{(n)}_{t}\>   \nonumber \\
  & & + \<\phi(\cdot) I_a^{-1}(\cdot), I_a(\cdot)Y^{(n)}_{t}\> + \<\phi(\cdot) I_a^{-1}(\cdot), I_a(\cdot) M^{(n)}_{t}\>,
 \eeqlb 
 for every $\phi \in {\cal S} (\R^d)$ and every $a\ge0$, 
 where, recalling that $\ell_\xi(s)=1_{ [\beta (\xi), \, \zeta (\xi))}(s)$,  \\
\beqnn
  & & U^{(n)}_{t}(\phi) :=  \frac{1}{\sqrt{{\theta}^{n}}}
\sum_{\xi \in {\Re}} \hspace{1mm}\sum_{p, i=1}^{d} \int_{0}^{t}
\ell_{\xi}(s)\,
\partial_p \phi(x_{\xi}(s))  c_{pi}(x_{\xi}(s))dB_{\xi i}(s)\,,  \\
&&   X^{(n)}_{t}(\phi) := \sum_{p=1}^{d}
\int_{0}^{t}\int_{\r^d}\<h_p(y-\cdot)
\partial_p \phi(\cdot), \mu^{(n)}_{s}\>W(dy,ds)\,, \\
\eeqnn
\beqnn
&&  Y^{(n)}_{t}(\phi) := \sum_{p,q=1}^{d} \int_{0}^{t}\left \<
\tfrac{1}{2} \partial_{p} \partial_{q} \phi(\cdot)\left[
\sum_{i=1}^{d}c_{pi}(\cdot)c_{qi}(\cdot) + \int_{\r^d}h_p(y -
\cdot)h_q(y - \cdot)dy \right], \mu^{(n)}_{s}\right\>ds,  \\
& &  M^{(n)}_{t}(\phi) :=  \int_{0}^{t} \int_{\r^d}
\phi(x)M^{(n)}(dx,ds) = \sum_{\xi \in {\Re}}
\frac{[O^{(n)}_{\xi}-1]}{{\theta}^{n}} \phi (x_{\xi}({\zeta(\xi)-}))
\, 1_{\{\zeta(\xi) \leq t\}} \,.
\eeqnn

The four terms in (\ref{DSPDE}) represent the respective contributions to the overall motion of the finite particle system
$\<\phi,\mu^{(n)}_{t}\>$ in $\R^d$ by the individual Brownian motions ($U^{(n)}_{t}(\phi)$), 
the random medium ($X^{(n)}_{t}(\phi)$), the mean effect of interactive and diffusive dynamics
($Y^{(n)}_{t}(\phi)$), and the branching mechanism ($M^{(n)}_{t}(\phi)$).  

Denote by $ \nu^{(n)}_t :=  T_{I_a}(\mu^{(n)}_{t})$ the transferred measure of interest and similarly define 
 \beqlb \lab{transf}
  v^{(n)}_t  :=  I_a(\cdot)U^{(n)}_{t} , x^{(n)}_t := I_a(\cdot)X^{(n)}_{t}, y^{(n)}_t := I_a(\cdot)Y^{(n)}_{t}, m^{(n)}_t :=  I_a(\cdot) M^{(n)}_{t}.
\eeqlb
 
We can now adapt Theorem 3.3 in \cite{RenSongWang09} to our present context. \\
 
 \bth \label{ApMain} 
  Assume Hypotheses $\ref{hyp:basicassumpFilter}$ and $\ref{hyp:basicassumpElliptic}$ are satisfied. 
  With the notation above, we have following conclusions:
\begin{description}

\item{\em (i)} $(\nu^{(n)}, v^{(n)}, x^{(n)}, y^{(n)}, m^{(n)})$ is tight on
$D([0, \infty), ({\cal S}^{\prime}(\R^d))^{5})$.

\item{\em (ii)} {\em (A Skorohod representation)}: Suppose that the joint distribution of a subsequence
\[
(\nu^{(n_m)}, v^{(n_m)}, x^{(n_m)}, y^{(n_m)}, m^{(n_m)}, W)
\]
converges weakly to the joint distribution of \[(\nu^{(0)}, v^{(0)}, x^{(0)}, y^{(0)}, m^{(0)}, W). \] 
Then, there exists a probability space $(\wt{\Omega}, \wt{\cal F}, \wt{\P})$ 
and $D([0, \infty),{\cal S}'(\R^d))$-valued sequences $\{\wt{\nu}^{(n_{m})} \}$,
$\{\wt{v}^{(n_{m})} \}$, $\{\wt{x}^{(n_{m})} \}$, $\{\wt{y}^{(n_{m})} \}$,
$\{\wt{m}^{(n_{m})} \}$ and a $D([0, \infty), {\cal S}' (\R^d))$-valued sequence $\{{\wt{W}}^{(n_{m})}\}$ defined on it,
such that
 \beqnn
&& \P \circ(\nu^{(n_m)}, v^{(n_m)}, x^{(n_m)}, y^{(n_m)}, m^{(n_m)}, W)^{-1} \\
&& \hspace{1cm}  = \wt{\P} \circ ( \wt{\nu}^{(n_{m})} ,
\wt{v}^{(n_{m})} , \wt{x}^{(n_{m})} , \wt{y}^{(n_{m})} , \wt{m}^{(n_{m})} ,
{\wt{W}}^{(n_{m})})^{-1}
 \eeqnn
 holds and, $\wt{\P}$-almost surely on $D([0, \infty),({\cal S}'(\R^d))^{6})$,
 \beqnn
\left( \wt{\nu}^{(n_{m})} , \wt{v}^{(n_{m})} , \wt{x}^{(n_{m})} , \wt{y}^{(n_{m})} ,
\wt{m}^{(n_{m})} , \wt{W}^{(n_{m})}\right) \rightarrow \left(
\wt{\nu}^{(0)} , \wt{v}^{(0)} , \wt{x}^{(0)} , \wt{y}^{(0)} , \wt{m}^{(0)} ,
\wt{W}^{(0)}\right) \quad
 \eeqnn
 as $ m \to \infty$. 
 
\item{\em (iii)} There exists a dense subset $\Xi \subset [0, \infty)$
such that $[0, \infty) \setminus \Xi$ is at most  countable  and for
each $t \in \Xi$ and each $\phi \in {\cal S}(\R^d)$, as $\R^{6}$-valued processes
 \beqnn
&& ( \wt{\nu}^{(n_{m})}_t(\phi) , \wt{v}^{(n_{m})}_t(\phi),
\wt{x}^{(n_{m})}_t(\phi), \wt{y}^{(n_{m})}_t(\phi), \wt{m}^{(n_{m})}_t(\phi),
{\wt{W}}^{(n_{m})}_t(\phi)) \\
&& \hspace{1cm}  \rightarrow ( \wt{\nu}^{(0)}_t(\phi),
\wt{v}^{(0)}_t(\phi), \wt{x}^{(0)}_t(\phi), \wt{y}^{(0)}_t(\phi), \wt{m}^{(0)}_t(\phi),
{\wt{W}}^{(0)}_t(\phi))
 \eeqnn
 in $L^{2}(\wt{\Omega}, \wt{\cal F}, \wt{\P})$ as $m \ra \infty$.
Furthermore, let $\wt{\cal F}^{(0)}_t$ be the $\sigma$-algebra
generated by $\{\wt{\nu}^{(0)}_{s}(\phi)$, $\wt{v}^{(0)}_{s}(\phi)$,
$\wt{x}^{(0)}_{s}(\phi)$, $\wt{y}^{(0)}_{s}(\phi)$, $\wt{m}^{(0)}_{s}(\phi)$,
${\wt{W}}^{(0)}_{s}(\phi)\}$ for all
$\phi \in {\cal S}(\R^d)$  and $ s \leq t$. Then
$\wt{m}^{(0)}_t(\phi)$ is a continuous, square-integrable $\wt{\cal
F}^{(0)}_t$-martingale with quadratic variation 
 \beqnn
\<\wt{m}^{(0)}_t(I^{-1}_a \phi) \> = \ga \sigma^2 \int_0^t \<I^{-1}_a \phi^2 ,
\wt{\nu}^{(0)}_u \>du.
 \eeqnn

\item{\em (iv)} $\wt{W}^{(0)}(dy,ds)$, ${\wt{W}}^{(n_{m})}(dy,ds)$
are Brownian sheets and, for any $\phi \in {\cal S}(\R^d)$, the continuous
square integrable martingale
\[
\displaystyle{ \wt{x}_{t}^{(n_{m})}(I^{-1}_a \phi):= \sum_{p=1}^{d}
\int_{0}^{t}\int_{\r^d} \left\<I^{-1}_a h_p(y-\cdot) \partial_{p}
\phi(\cdot), \wt{\nu}^{(n_{m})}_{s}\right\>
{\wt{W}}^{(n_{m})}(dy,ds) }
\]
converges to
$$
\wt{x}_{t}^{(0)}(I^{-1}_a \phi)  := \sum_{p=1}^{d} \int_{0}^{t}\int_{\r^d}
\left\< I^{-1}_a h_p(y-\cdot) \partial_p \phi(\cdot), \wt{\nu}^{(0)}_{s}
\right\> \wt{W}^{(0)}(dy,ds) $$ in $ L^{2}(\wt{\Omega}, \wt{\cal F},
\wt{\P})$ .

\item{\em (v)} $\wt{\nu}^{(0)}=\{ \wt \nu^{(0)}_t: t\geq 0\}$ is
a solution to the $({\cal L},\delta_{\nu_0})$-martingale problem and
$\wt{\nu}^{(0)}$ is a continuous process and for every $t \ge 0$ and every $\phi \in
{\cal S}(\R^d)$ we have
 \beqlb \label{SPDE}
 & & \wt{\nu}^{(0)}_{t}(I^{-1}_a \phi) - \wt{\nu}^{(0)}_{0}(I^{-1}_a \phi) =
\wt{x}^{(0)}_{t}(I^{-1}_a \phi) \nonumber \\
& &  + \int_{0}^{t} \left\< \sum_{p,q=1}^{d}
\tfrac{1}{2}I^{-1}_a (a_{pq}(\cdot)+\rho_{pq}(\cdot, \cdot)) \partial_p
\partial_q \phi(\cdot), \, \wt{\nu}^{(0)}_{s} \right\>ds   \nonumber \\
& &\hspace{1cm} + \int_{0}^{t}\int_{\R^d}I^{-1}_a \phi(x)
\wt{m}^{(0)}(dx,ds)
\mbox{\hspace{1cm}$\wt{\P}$-a.s. and  in $L^{2}(\wt{\Omega}, \wt{\cal F}, \wt{\P})$}.
 \eeqlb
\end{description}
 \eth
 \proof 
 While this result is not a direct consequence of Theorem 3.3 in \cite{RenSongWang09}, our conditions ensure 
 that all the key inequalities used in their proof remain valid here as well. Some of the difficulties encountered there,  
 due to the necessity of dealing with trajectories in a much larger space, are here avoided through the direct use of 
 ${\cal S}' (\R^d)$, the dual of a nuclear Fr\'{e}chet space. 
 We therefore only summarize those key ideas that fit our own context. 
 
(i) By Mitoma \cite{Mitoma83}, we only need to prove that, for any $\phi \in {\cal S}(\R^d)$,  
the sequence of laws of $(\nu^{(n)}(\phi), v^{(n)}(\phi), x^{(n)}(\phi),  y^{(n)}(\phi), m^{(n)}(\phi))$ 
is tight in $D([0, \infty), \R^{5})$, after noting that $\phi\cdot I_a^{-1}\in{\cal S}(\R^d)$ holds. 
The arguments in the proof of Lemma 3.2 in \cite{RenSongWang09} still work and we can use the tightness criterion 
from Theorem 3.8.6 in Ethier-Kurtz \cite{EthierKurtz86} to obtain the tightness for each $\phi$.

(ii) This is a direct application of the versions of Skorohod's Representation Theorem 
on $D([0, \infty), ({\cal S}^{\prime}(\R^d))^k)$ due to Jakubowski \cite{Jakubowski86,Jakubowski97}. 
See also \cite{CWX12} in the present context. 

(iii) The proof is the same as that of part (iii) of Theorem 3.3 in \cite{RenSongWang09}.
 
(iv) Since $W$, ${\wt{W}}^{(0)}$ and ${\wt{W}}^{(n_{m})}$ have the
 same distribution, $\wt{W}^{(0)}$ and ${\wt{W}}^{(n_{m})}$ are
 Brownian sheets. The conclusion follows from (ii) just as in \cite{RenSongWang09}. 
 
(v) An application of It\^{o}'s formula under conditions that remain unchanged from the corresponding case treated 
in \cite{RenSongWang09} yield that $\{\wt{\mu}_{t}^{(0)}: t \geq 0\}$ is a solution to the martingale problem for 
$({\cal{L}},\delta_{\mu_{0}})$ with continuous trajectories.
\qed 

This proves part of Theorem \ref{MPforSDSM}, namely the existence of a solution to the $({\cal L}, \delta_{{\mu}_{0}})$-martingale problem
for the operator given by $(\ref{pregenerator})$, on $C([0, \infty), {\cal S}'(\R^d))$
for any starting measure $\mu_0\in M_a(\R^d)$. It also confirms the validity of both $(\ref{SPDEc})$ and $(\ref{SPDEb})$, 
as well as the fact that both ${X}_{t}(\phi)$ and $M_t(\phi)$
are mutually orthogonal continuous square-integrable martingales. 
It remains to show that there is one such solution living on $C([0, \infty), M_a(\R^d))$ and that this solution is unique. 
This is achieved next. 


\section{Duality and the martingale problem}\label{sec:dualConst}
\setcounter{equation}{0}
In this section, we use the construction of a function-valued dual process, in the sense of Dawson and Kurtz \cite{DawsonKurtz82},
as a way to directly exhibit the transition probability of SDSM, thus immediately giving an alternative construction of SDSM, 
as well as proving uniqueness in law on the space $C([0, \infty), M_a(\R^d))$, since duality also yields the full characterization of the law of SDSM
by way of the martingale problem formulation. The technique of duality was developed in order to identify the more complex measure-valued one uniquely and compute some of its mathematical features, after first showing its existence through some other means, often by way of a tightness argument or some other limiting scheme. Part of the interest of this section lies with the uncommon use of the existence of a dual function-valued process, in order to construct a transition function for SDSM and show the existence of associated measure-valued processes of interest on richer spaces of trajectories resulting from the inclusion of infinite starting measures. 
The approach used here was developed by Barton et al. \cite{BartonEtheridgeVeber10} for the construction and characterization 
of a spatial version of the $\Lambda$-Fleming-Viot process, after the work of Evans \cite{Evans97} on systems 
of coalescing Borel right processes. See also Donnelly et al. \cite{DEFKZ2000} for another instance of this, 
in the context of the continuum-sites stepping-stone process. 
(Note that some technical aspects of the treatment in this section are required
due to the topology on $M_a(\R^d)$. The reader can refer to the appendix in
Konno and Shiga \cite{KonnoShiga88} for additional clarifications.) 

Let us begin with the construction of the function-valued process that will serve our purpose,
namely an extension of the ones
built in Ren et al. \cite{RenSongWang09} and Dawson et al. \cite{DLW01}. 

In order to facilitate some of the calculations required henceforth,
notably because infinite starting measures lying in $M_a(\R^d)$ impose
restrictions on the set of functions needed for a full description of the dual process,
the domain ${\cal D}({\cal L})$ of operator ${\cal L}$ in (\ref{pregenerator})
 --- the set of functions in $B(M_a(\R^d))$ upon which ${\cal L}$ is well-defined ---
is enlarged to comprise all bounded continuous functions of the form
 \beqlb\label{eqn:fulldomain}
F(\mu) = g(\< f_1, \mu^{m_1} \>,  \cdots, \< f_k, \mu^{m_k} \>)
\eeqlb
with $g \in C^2(\R^k)$ for some $k\ge1$,
any choice of positive integers $m_1,\ldots, m_k$ and, for every $1 \leq i \leq k$,
$f_i  \in {\cal D}_a(G_{m_i})$. 

We describe the space ${\cal D}_a(G_{m})$ next.
For the generator $G_m$ from (\ref{eqn:Gn}) of strongly continuous contraction semigroup
$\{P^{m}_t\}$ on Banach space $C_0((\R^d)^{m})$, the domain ${\cal D}(G_m)$
--- the set of functions in $B((\R^d)^{m})$ upon which $G_m$ is well-defined ---
is simply the set of those functions $f$ such that the limit
\[
\lim_{t \ra 0+}\frac{1}{t}(P^{m}_t f - f)
\]
exists in $C_0((\R^d)^{m})$, so we write $f \in {\cal D}(G_{m})$ if and only if this limit exists and equals $G_{m}f$. \\
In order to ensure integrability with respect to some infinite measures,
our statements about functions in this domain ${\cal D}(G_m)$ are restricted to its subspace
defined by
\beqlb \label{DaGm}
 {\cal D}_a(G_m):= \{ f \in {\cal D}(G_m) : \| {\cal I}_{a, m}^{-1} f \|_{\infty} < \infty
 \mbox{ and  $\| {\cal I}_{a, m}^{-1} G_m f \|_{\infty} < \infty  $} .\}
 \eeqlb

The short form $\mu^{m}= \mu\otimes\ldots\otimes\mu$
denotes the $m$-fold product measure of $\mu\in M_a(\R^d)$ by itself
and we write ${\cal I}_{a, m}$ for the product
${\cal I}_{a, m}(x)=I_a(x_1)\cdot\ldots\cdot I_a(x_m)$,
keeping in mind that
${\cal I}_{a, m}^{-1} f(x)={\cal I}_{a, m}^{-1}(x)f(x)$
means the product, not the composition of functions. \\

Observe first that both
 \beqlb\label{iam}
{\cal I}_{a, m} \in C_b^{\infty}((\R^d)^{m})
\eeqlb
and
 \beqlb\label{iamgm}
{\cal I}_{a, m}^{-1} G_m {\cal I}_{a, m} \in C_b((\R^d)^{m})
\eeqlb
hold, under Hypothesis $\ref{hyp:basicassumpElliptic}$,
hence so do $\| {\cal I}_{a, m}^{-1} G_m {\cal I}_{a, m} \|_{\infty} < \infty$
and ${\cal I}_{a, m} \in {\cal D}_a(G_m)$.
A quick sketch of proof of these facts is supplied in Subsection \ref{app:pf_iam}.
The useful inclusions $C_c^2((\R^d)^m) \subset{\cal D}_a(G_m)\subset {\cal D}(G_m)$
and ${\cal I}_{a, m}^{-1} G_m \{C_c^2((\R^d)^m)\} \subset C_c((\R^d)^m)$
are also clearly valid for every choice of $a\ge0$.

It is important to note at this point that, for every positive value of $a > 0$
 and $m\ge1$, while ${\cal I}_{a, m} \in C_0^{\infty}((\R^d)^{m})$ holds (this is false when $a=0$),
 we also have ${\cal I}_{a, m} \not\in {\cal D}_b(G_m)$ for any $b > a$. Therefore
 $C_0^{\infty}((\R^d)^m) \not\subset {\cal D}_a(G_m)$ for any $a > 0$,
so the core $C_0^2((\R^d)^m)$ of $G_m$ does not lie inside ${\cal D}_a(G_m)$
even though $C_c^2((\R^d)^m)$ is uniformly dense in $C_0^2((\R^d)^m)$.

More generally, we also get the following results
pertaining to the preservation of the semigroup property under the rescaling
induced by function ${\cal I}_{a, m}$, the proof of which may also be found in Subsection \ref{app:pf_iam}.

 \ble \label{lea}
 Assume that Hypothesis $\ref{hyp:basicassumpElliptic}$ is satisfied.
 For every $a\ge0$, $f \in {\cal D}_a(G_m)$ and $T > 0$, there holds 
 $P_T^{m} f \in  {\cal D}_a(G_m)$,
$\sup_{0 \leq t \leq T} \| {\cal I}_{a, m}^{-1} P_t^{m} f \|_{\infty} < \infty$ and
\[
\sup_{0 \leq t \leq T} \| {\cal I}_{a, m}^{-1} \frac{\partial}{\partial t} P_t^{m} f \|_{\infty}
=\sup_{0 \leq t \leq T} \| {\cal I}_{a, m}^{-1} G_m P_t^{m} f \|_{\infty}
=\sup_{0 \leq t \leq T} \| {\cal I}_{a, m}^{-1} P_t^{m} G_m f \|_{\infty} < \infty.
\]
 \ele

The construction of the function-valued process can now proceed, as follows.

Let $\{J_t: t\ge 0\}$ be a decreasing c\`{a}dl\`{a}g Markov jump process
on the nonnegative integers $\{0,1,2,\ldots\}$, started at $J_0=m$
and decreasing by 1 at a time,
with Poisson waiting times of intensity $\gamma \sigma^2 l(l-1)/2$
when the process equals $l\ge2$. The process is frozen in place
when it reaches value $1$ and never moves if it is started at
either $m=0$ or $1$.
Write $\{\tau_k: 0\le k\le J_0-1\}$ for the sequence of jump times of
$\{J_t:t\ge 0\}$ with $\tau_0=0$ and $\tau_{J_0} = \infty$.
At each such jump time a randomly chosen projection is effected on the
function-valued process of interest, as follows.
Let $\{S_k: 1\le k\le J_0\}$ be a sequence of random
operators which are conditionally independent given $\{J_t: t\ge 0\}$ and satisfy
$$
\P\{S_k = \Phi_{ij}^m | J_{\tau_k -} =m\} = \frac{1}{m(m-1)}, \qquad
1 \le i \neq j \le m,
$$
as long as $m\ge2$.
Here $\Phi_{ij}^mf$ is a mapping from ${\cal D}_a(G_m)$ into ${\cal D}_a(G_{m-1})$ defined by
\beqlb \label{restriction}
\Phi_{ij}^m f(y):= f(y_{1}, \cdots, y_{j-1},y_{i},y_{j+1},\cdots, y_{m}),
\eeqlb
for any $m\ge2$ and $y=(y_1, \cdots, y_{j-1}, y_{j+1}, \cdots, y_m)\in (\R^d)^{m-1}$.
When $m=0$ we simply write ${\cal D}_a(G_0)=\R^d$ and $P_t^0$ acts as
the identity mapping on constant functions.

That $\Phi_{ij}^m$ is well-defined follows from the observation that the sets ${\cal D}_a(G_m)$ form
an increasing sequence in $m$, in this last case when interpreting any
function of $m\le n$ variables also as the restriction of a function of $n$ variables.
Details are in Subsection \ref{app:pf_iam}.

Given $J_0=m$ for some $m\ge0$, define process
$Y:=\{Y_t: t\ge0\}$, started at some point $Y_0\in {\cal D}_a(G_m)$ within
the (disjoint) topological union $\bB:=\cup_{m=0}^\infty {\cal D}_a(G_m)$, by
\beqlb \label{Yprocess}
Y_t = P^{J_{\tau_k}}_{t-\tau_k} S_k P^{J_{\tau_{k-1}}} _{\tau_k
-\tau_{k-1}}S_{k-1} \cdots P^{J_{\tau_1}}_{\tau_2 -\tau_1} S_1
P^{J_0}_{\tau_1}Y_0, \quad \tau_k \le t < \tau_{k+1}, 0\le k\le J_0-1.
\eeqlb
By Lemma \ref{lea}, the process $Y$ is a well-defined $\bB$-valued
strong Markov process for any starting point $Y_0\in \bB$.
Clearly, $\{(J_t, Y_t): t\ge 0\}$ is also a strong Markov process.

\ble \label{lea2}
Assume that Hypothesis $\ref{hyp:basicassumpElliptic}$ is satisfied.
Given any $a\ge0$, $J_0=m\ge1$ and $T > 0$,
there exists a constant $c=c(a,d,m,T) > 0$ such that,
for every $Y_0\in {\cal D}_a(G_m)$ we have $\P$-almost surely
\[
 \sup_{0 \leq t \leq T} \| {\cal I}_{a, J_t}^{-1} Y_t \|_\infty \le c \| {\cal I}_{a, m}^{-1} Y_0 \|_\infty.
\]
\ele

The proof is found in Subsection \ref{app:pf_iam}. 

Recall from (\ref{transfer}) the homeomorphic transfer function $T_{I_a} : M_a(\R^d) \rightarrow M_0(\R^d)$ denoted by
$$ 
T_{I_a}(\mu)(A) := \int_{A} I_a(x) \mu(dx)=\int_{A} (1 + |x|^2)^{- a/2} \mu(dx),
$$
for any $A \in {\cal B}(\R^d)$.\\

\bth \label{tha}
  Assume that Hypotheses $\ref{hyp:basicassumpFilter}$ and $\ref{hyp:basicassumpElliptic}$ are satisfied. For any $a\ge0$, $m\ge1$,
  $ f\in {\cal D}_a(G_m)$, $\mu_0\in M_a(\R^d)$ and $t\in[0,\infty)$, there exists a
  time-homogeneous Markov transition function
  $\{{\Q}_t(\mu,\Gamma): t\in[0,\infty),\mu\in M_a(\R^d),\Gamma\in {\B}(M_a(\R^d))\}$,
 given by
\beqlb \label{dual}
 & &  \int_{M_a(\R^d)}\<f, \nu^m\>{\Q}_t(\mu, d\nu)  \nonumber \\
 = &   & \E \left[\<{\cal I}_{a, J_t}^{-1} Y_t, (T_{I_a}(\mu))^{J_t}\>
\exp\left(\frac{\gamma \sigma^2}{2}\ \int_0^t  J_s(J_s-1)ds\right)\Bigl|(J_0,Y_0)=(m,f) \right].
\eeqlb
The associated probability measure $\P_{\mu_0}$ on $C([0, \infty), M_a(\R^d))$ of the form:
 \beqlb \lab{GSDSM}
  &&
 \P_{\mu_0} (\{ w \in C([0, \infty), M_a(\R^d)): w_{t_i} \in \Gamma_i, i=0, \cdots, n\})          \nonumber   \\
 = && \int_{\Gamma_0} \cdots \int_{\Gamma_{n-1}}
   {\Q}_{t_n-t_{n-1}}(\mu_{n-1}, \Gamma_n){\Q}_{t_{n-1}-t_{n-2}}(\mu_{n-2}, d \mu_{n-1}) \cdots {\Q}_{t_{1}}(\mu_{0}, d \mu_{1}),
 \eeqlb
 for any $0 \leq t_1 \leq t_2 \leq \cdots \leq t_n$ and $\Gamma_i \in {\cal B}(M_a(\R^d)), i=0,1, \cdots, n$, 
 is the unique probability measure on $C([0, \infty), M_a(\R^d))$ which satisfies the 
 (well-posed) $({\cal L}, \delta_{{\mu}_{0}})$-martingale problem.
\eth
 \noindent  {\bf Remark:}
 We adapt the approach used for the case $a=0$ and $d=1$ in Dawson et al. \cite{DLW01} by exhibiting
 a transition function, built by using the law of function-valued process $Y$ and charging
 space $C([0, \infty), M_a(\R^d))$ with a probability measure fitting our needs.
 In this setting, we only give a quick sketch of the main ideas but provide details for
 overcoming the new difficulties arising from the larger space.

 \proof
 We need to prove that the time-homogeneous transition function given by (\ref{dual}) is well-defined, so that
 ${\Q}_t(\mu, \cdot)$ is a probability measure on ${\B}(M_a(\R^d))$ with ${\Q}_0(\mu, \cdot)=\delta_{\mu}$
 and ${\Q}_{\cdot}(\cdot,\Gamma)$ is Borel measurable, for every choice of $t$ and $\mu$.

 We begin by quickly sketching that there exists
 a transition function ${\Q}_t(\mu, d \nu)$ on ${\B}(M_0(\R^d))$, for any $t \geq 0$ and $\mu \in M_0(\R^d)$,
 for which the duality identity (\ref{dual}) and the extended Chapman-Kolmogorov equations (\ref{GSDSM}) both hold.
 The argument for this first step is as follows.

 In the general case $a\ge0$, Lemma \ref{lea2} ensures that
 the right hand side of (\ref{dual}) is well-defined, as well as bounded for every $ f\in {\cal D}_a(G_m)$;
 therefore, it defines a strongly continuous positive bounded semigroup of linear operators
 on the subspace of $C_b(M_a({\R^d}))$ of functions of the form
 $F_{m,f} (\mu):=\< f, \mu^m \>$ for any $f \in {\cal D}_a(G_m)$. 
 (The argument runs along the same lines as in the proof of Theorem 5.1 in Dawson et al. \cite{DLW01} for the case $d=1$.)
 Unfortunately the constant functions, ${\cal D}_a(G_0)$, do not have finite integrals when $a > 0$
 and must therefore be removed from the set  $\bB:=\cup_{m=0}^\infty {\cal D}_a(G_m)$ of values,
so this amputated subspace of functions suffices no longer for a full characterization of measures.

Since $M_a({\R^d})$ is not locally compact, one needs, at this stage,
 to first add a point $\infty$ at infinity to build the one-point compactification
 $\hat{\R}^d:=\R^d\cup\{\infty\}$ of $\R^d$ and then add an isolated point $\Delta$
 to the space of finite Radon measures $M_0(\hat{\R}^d)$ on $\hat{\R}^d$,
 yielding the new compact metrizable space $E^{\Delta}:=M_0(\hat{\R}^d)\cup\{\Delta\}$.
 By the Stone-Weierstrass theorem, the algebra $\Pi(E^{\Delta})$ generated by the set of all functions of the form
 $F_{m,f} (\mu):=\< f, \mu^m \>$, for any $m\ge1$ and $f \in C((\hat{\R}^d)^m)$, is dense in $C(E^{\Delta})$.
 The operator ${\cal T}_tF_{m,f}$ defined as the right hand side of (\ref{dual}) can be extended to $C(E^{\Delta})$
 by simply putting ${\cal T}_tF:=F(\Delta)+{\cal T}_t(F-F(\Delta))$ for every $F\in\Pi(E^{\Delta})$.

 Theorem 4.2.7 in Ethier and Kurtz \cite{EthierKurtz86} implies
 the existence of a probability measure ${\Q}_t(\mu, \cdot)$ on ${\B}(E^{\Delta})$ such that (\ref{dual}) holds,
 for every $f \in C((\hat{\R}^d)^m)$, $\mu_0\in E^{\Delta}$, $\mu\in E^{\Delta}$, $\Gamma\in{\B}(E^{\Delta})$ and $t\in[0,\infty)$. 
 The definition of the extension ensures that ${\Q}_t$ is normalized so that it is a probability. 
 The existence of a probability measure $\P_{\mu_0}$ on $D([0, \infty), E^{\Delta})$ which satisfies (\ref{GSDSM}) 
 also ensues from Theorem 4.2.7 in Ethier and Kurtz \cite{EthierKurtz86}. That $\P_{\mu_0}$ is supported by 
 $C([0, \infty), E^{\Delta})$ follows from the results of Bakry and Emery \cite{BakryEmery85} on operators with the derivation proprety, 
 just as in the special case $d=1$ treated in Wang \cite{Wang98}. 
 
 The canonical process $\{w_{t}\}$ under $\P_{\mu_0}$ is our SDSM
 with initial distribution $\delta_{\mu_0}$ and $\Delta$ as a trap.
 Just as in the original case $a=0$ and $d=1$ treated in Wang \cite{Wang98} (see his Theorem 4.1),
 unicity of the solution for (\ref{GSDSM}) on $C([0, \infty), E^{\Delta})$, started at $\mu_0\in E^{\Delta}$,
 follows from the well-posedness of the $({\cal L}, \delta_{{\mu}_{0}})$-martingale problem,
 by way of the duality identity (\ref{dual}) on $C([0, \infty), E^{\Delta})$
plus the fact that, for any $t > 0$, the moment power series
 \[
 \sum_{n=1}^{\infty}\frac{\theta^n}{n!} \E_{\mu_0} \< 1, w_t\>^n
 \]
has a positive radius of convergence, thus identifying the one dimensional laws uniquely.
This works because we chose to build $E^{\Delta}$ with $a=0$ rather than $a > 0$.

By an argument similar to
Step 8 in the appendix of Konno-Shiga \cite{KonnoShiga88} or the proof of Theorem 4.1 of Dawson et al \cite{DLW01},
we get $\P_{\mu_0}(w_{t}(\{\infty\}))=0 \mbox{    for all $t \geq 0$})= 1$ provided ${\mu_0}(\{\infty\})=0$.
It follows that whenever $\mu_0\in M_0(\R^d)$ holds instead of $\mu_0\in E^{\Delta}$, the solution
$\P_{\mu_0}$ is actually supported by $C([0, \infty), M_0(\R^d))$.

This ends the proof of both (\ref{dual}) and (\ref{GSDSM}) in the case $a=0$.\\

With $a > 0$, given $\{{\Q}_t(\mu,\Gamma): t\in[0,\infty),\mu\in M_0(\R^d),\Gamma\in {\B}(M_0(\R^d))\}$
just constructed, define
$\{{\tilde {\Q}}_t(\mu,\Gamma): t\in[0,\infty),\mu\in M_a(\R^d),\Gamma\in {\B}(M_a(\R^d))\}$ by
 \[
{\tilde {\Q}}_t(\mu,\Gamma):={\Q}_t(T_{I_a}(\mu),T_{I_a}(\Gamma))
 \]
and ${\tilde \P}_{\mu_0}$ on $C([0, \infty), M_a(\R^d))$ by (\ref{GSDSM})
with ${\tilde {\Q}}_t$ instead of ${\Q}_t$ throughout and a typical trajectory written as
$\{{\tilde w}_t\} \in C([0, \infty), M_a(\R^d))$.
Both ${\tilde {\Q}}_t$ and ${\tilde \P}_{\mu_0}$ are well-defined solutions to (\ref{dual}) and (\ref{GSDSM}) respectively.
Since $M_a(\R^d)$ is a Polish space, the uniqueness of solution to the martingale problem 
is equivalent to that of $\{{\Q}_t:t\ge0\}$, by Theorem 4.1.1 of Ethier and Kurtz \cite{EthierKurtz86}.
Both hold because the homeomorphism $T_{I_a}$ allows us to write,
for any two solutions $\{{\tilde {\Q}}_{\xi,t}:t\ge0\}$ to (\ref{dual}), for $\xi=1,2$,
with $\{{\Q}_t:t\ge0\}$ the unique solution in the case $a=0$:
 \[
 \int_{M_a(\R^d)}\<f, \nu^m\>{\tilde {\Q}}_{\xi,t}(\mu_0, d\nu)
 = \int_{M_0(\R^d)}\<{\cal I}_{a, m}^{-1}f, T_{I_a}(\nu)^m\> {\Q}_t(T_{I_a}(\mu_0), d[T_{I_a}(\nu)])
\]
for any $m\ge1$, $ f\in {\cal D}_a(G_m)$, $\mu_0\in M_a(\R^d)$ and $t\in[0,\infty)$.
Since $K_a(\R^d)\subset {\cal D}_a(G_1)$ is a separating set
 for measures inside $M_a(\R^d)$, it ensues that $\cup_{m=1}^\infty {\cal D}_a(G_m)$
 is a algebra of functions that separates points in $M_a(\R^d)$ and hence is separating
 for the probability measures on ${\cal B}(M_a(\R^d))$,
 by Theorem 3.4.5 of Ethier and Kurtz \cite{EthierKurtz86}.

Finally, denoting by ${\cal L}^{*}$ the infinitesimal generator for the process $(J,Y)$,
the differential form of the duality between its associated martingale problem and
the $({\cal L}, \delta_{{\mu}_{0}})$-martingale problem of pregenerator (\ref{pregenerator}) for SDSM
can be written in the following manner. We only need to consider
functions of the form (\ref{eqn:fulldomain}) with $k=1$, since they form a separating set
for measures in $M_a(\R^d)$, by Proposition 4.4.7 in Ethier and Kurtz \cite{EthierKurtz86}.
Writing $F_{m,f} (\mu) := g(\< f, \mu^m \>)$ for any such function with $g(s)=s$,
$f \in {\cal D}_a(G_m)$ and $\mu \in M_a(\R^d)$,
a quick calculation yields
\beqlb \label{dualityform}
{\cal L} F_{m,f} (\mu)=\< {\cal L}^{*}f, \mu^m \> + \frac{1}{2} \ga \sigma^2 m(m-1)\< f, \mu^m \>
\eeqlb
with
\[
{\cal L}^{*}f={G}_{m} f + \frac{\gamma \sigma^2}{2} \sum_{i,j =1, \, i\not=j}^m (\Phi_{ij}^mf - f).
\]
An application of  It\^{o}'s formula with a stopping argument allows this formal calculation
to be made rigorous, so that (\ref{SPDEb}) ensues for every $\phi\in C_c^\infty(\R^d)$.
See Ren et al. \cite{RenSongWang09} for additional details.
We need only check (\ref{SPDEb}) for the remaining function $\phi=I_a$ to cover all of $K_a(\R^d)$.
Since $I_a\in L^1(\mu_0)$ holds by definition,  $I_a^{-1} G_1 I_a$ satisfies (\ref{iamgm})
and $\partial_p I_a\in C_0(\R^d)\cap L^1(\mu_0)$
by the calculations in Subsection \ref{app:pf_iam}, the same stopping argument can be applied,
as the three integrals in (\ref{SPDEb}) can be controlled simultaneously. 

Therefore, the uniqueness follows by applying the now classical results from Dawson and Kurtz \cite{DawsonKurtz82}
(or use Corollary 4.4.13 in Ethier and Kurtz \cite{EthierKurtz86}) to the generators appearing in the
differential form (\ref{dualityform}) of the duality equation, the required integrability conditions holding 
for every $\phi\in K_a(\R^d)$, which is uniformly dense in ${\cal D}_a(G_1)$. 
\qed

\vfill
\break

 \proof [Proof of Theorem \ref{MPforSDSM}] 
Theorem \ref{tha} yields the existence and uniqueness of the solution to the martingale problem for SDSM 
on the space $C([0, \infty), M_a(\R^d))$. 
A representation of this unique solution has already been obtained in Theorem \ref{ApMain} in the form of the 
solution to an SPDE on $C([0, \infty), {\cal S}'(\R^d))$. 
This representation satisfies all of the other statements in Theorem \ref{MPforSDSM} through the continuity of 
the transfer function $T_{I_a}$ defined in (\ref{transfer}) and that of the injection of 
$M_0(\R^d)$ in ${\cal S}' (\R^d)$, by Proposition 5.1 in Dawson and Vaillancourt \cite{DV95}.
\qed 

Another consequence of Theorem \ref{tha} is the following uniform bound on the growth of the moments associated with 
the single particle transition density $q_t^{1}(\cdot,x)$ with generator $G_1$ on $\R^d$ given by $(\ref{eqn:Gn1})$ and its spatial derivative 
$\partial_p q_t^{1}(\cdot,x)$. 

   \bcor \lab{MCOM} Assume all the conditions in of Theorem \ref{tha}. Let $\{\mu_t\}$ be the SDSM started at a $\mu_0$ 
   which satisfies Hypothesis $\ref{hyp:basicassumpGauss}$, with dual process $\{Y_t\}$ 
   started at $Y_0=f$. The duality identity (\ref{dual}) holds true for all functions $f\in L^1(\mu_0^m)$,
with a finite common value on both sides as part of the conclusion; these include the two cases : 
a) $f=q_\epsilon^{1}(\cdot,x)$ and b) $f=\partial_p q_\epsilon^{1}(\cdot,x)$, for any fixed $x\in\R^d$ and $\epsilon>0$.
  The dual identity (\ref{dual}) can also be rewritten as follows : for every integer $m>0$,
  \beqlb \lab{ME}
  \E_{\mu_0}[\<f, \mu_t\>^m] = \E_{m,f^{\otimes m}} \bigg[\<Y_t, \mu_0^{J_t}\> \exp\{\frac{1}{2} \int_0^t J_s(J_s-1)ds \}  \bigg]
  = \sum_{k=0}^{m-1}M_k(t),
  \eeqlb
with $\sup_{0\le t\le T}|M_k(t)|<\infty$ for every $T>0$, where process $Y$ is expressed as in (\ref{Yprocess}) and 
  \beqlb \lab{mMk0}
 M_k(t)  & :=  & \frac{m! (m-1)!}{2^k (m-k)! (m-k-1)!} \int_{(0, t]}dr_1 \int_{(r_1, t]}dr_2 \cdots \int_{(r_{(k-1)}, t]}  \nonumber \\
      &     & \cdot  \E_{m,f^{\bigotimes m}}\bigg[ \<P^{m-k}_{t-r_k} S_k \cdots  P^{m-1}_{r_2-r_1} 
      S_1  P^{m}_{r_1}f^{\bigotimes (m)}, \mu_0^{\bigotimes (m-k)}\>\bigg| \tau_j = r_j : 1 \leq j \leq k \bigg] dr_k \nonumber \\ 
  \eeqlb
  \ecor
  
  \proof 
  First, recall that all Radon measures are locally finite, 
  hence all compact sets have finite measure and 
  $C_c^\infty(\R^d)\subset L^p(\mu_0)$, for all $a\ge0$ and $p\ge1$. 
  Since all Radon measures are inner regular and outer regular, for any set $N\in{\cal B}(\R^d)$ we can choose a decreasing 
 sequence of open subsets $U_\ell\in{\cal B}(\R^d)$ and an increasing sequence of compact 
 subsets $V_\ell\in{\cal B}(\R^d)$ such that $V_\ell \subset N \subset U_\ell$
 and $\mu_0(U_\ell\smallsetminus V_\ell) < 1/\ell$. In this context, 
 Urysohn's Lemma (see Lemma 2.19 in Lieb and Loss \cite{LiebLoss01}) asserts the
 existence of a sequence of functions $\psi_\ell\in C_c^\infty(\R^d)$ such that 
 $1_{V_{\ell}}\le\psi_\ell\le 1_{U_\ell}$ holds everywhere, with $1_N(x)=1$ if $x\in N$ and $0$ elsewhere. 
 This proves not only that $C_c^\infty(\R^d)$ is dense in $L^p(\mu_0)$, for all $p\ge1$, 
 but also that every $\phi\in L^p(\mu_0)$ is the $\mu_0$-almost sure limit as well as the 
 $\|\cdot\|_{\mu_0,p}$-norm limit of some sequence $\phi_n\in C_c^\infty(\R^d)$. 
 (By the way, note that we have just proved that $C_c^\infty(\R^d)$ is  $\P_{\mu_0}$-almost surely dense 
 in $L^p(\mu_t)$, for all $t > 0$, $a\ge0$ and $p\ge1$ as well,
 since Theorem \ref{MPforSDSM} implicitly asserts that $\mu_t\in M(\R^d)$ holds $\P_{\mu_0}$-almost surely.) 
 Meanwhile, $Y_t\in C_0^2((\R^d)^m)$ holds for all $t>0$ as soon as $Y_0\in C_c^\infty((\R^d)^m)$ does.
 Extending the duality identity (\ref{dual}) through a density argument therefore makes sense.
 We do the proof only for $\mu_0\in M_0(\R^d)$ and leave it to the reader to extend it to $\mu_0\in M_a(\R^d)$ for $a>0$.

 The duality identity (\ref{dual}) first extends to all $f\in B((\R^d)^m)\cap L^1(\mu_0^m)$
since $C_0^2((\R^d)^m)$ is uniformly dense in $C_0((\R^d)^m)$ and therefore dense in 
$B((\R^d)^m)$ in the bounded pointwise sense 
(see Proposition 3.4.2 in Ethier and Kurtz \cite{EthierKurtz86}). 
Moreover, (\ref{dual}) remains true for unbounded $f\in L^1(\mu_0^m)$, by monotone convergence
of the sequence $f\cdot 1_{\{| f | < n\}}$ applied to the positive and negative parts of $f$ on both sides 
of (\ref{dual}), since operators $P_t^{m}$ and $\Phi_{ij}^m$ are both linear and positivity preserving. 
When $r>0$ and $x\in\R^d$, $f \in L^1(\mu_0)$ in both the cases considered here 
--- $f=q_\epsilon^{1}(\cdot,x)$ and $f=\partial_p q_\epsilon^{1}(\cdot,x)$ --- by Hypothesis $\ref{hyp:basicassumpGauss}$ and $(\ref{LSU})$, 
so the integrals on both sides of the first equality in $(\ref{ME})$ 
make sense in those cases as well. Conditionning on the number of jumps for dual process $\{Y_t\}$ yields $(\ref{mMk0})$. 
It remains to show $\sup_{0\le t\le T}|M_k(t)|<\infty$ and that the expectations involved in the above argumentation, are all finite.
This is done on the right hand side of (\ref{dual}) by reducing the proof for the SDSM case, to that of Super-Brownian motion. 
   
   First note that all four constants $a^*$, $b$, $c$ and $A^*$ in the Aronson bounds $(\ref{Aronsonbounds})$ 
   depend on the dimension $d$ of the space. Without loss of generality we can select them so that $(\ref{Aronsonbounds})$ 
   holds jointly for all $d\le m$ with common values for $a^*$, $b$, $c$ and $A^*$ --- the proof is found in Subsection \ref{app:pf_jointparametrization}. 
   To keep notation to a minimum in this proof, let $K>0$ denote some generic constant, which will change throughout the proof but the 
   value of which is irrelevant for the sake of this argumentation --- its value depends on $a^*$, $b$, $c$ and $A^*$.  
   
   Writing $T^k_t$ for $P^k_t$ in the special case of the heat semigroup and using the exchangeability of the dual particles in this particular instance,   
    $(\ref{Aronsonbounds})$ implies, for any nonnegative $f \in L^1(\mu_0)$ and any choice of $0<r_1<r_2<\ldots<r_k<t\le T$,  
 \beqlb \lab{mMk}
 &&  \<P^{m-k}_{t-r_k} S_k \cdots  P^{m-1}_{r_2-r_1} S_1  P^{m}_{r_1}f^{\bigotimes (m)}, \mu_0^{\bigotimes (m-k)}\>  \nonumber \\
  \le && K  \< ( T^1_{c t  -  c r_k } )^{\bigotimes (m-k)} S_k \cdots  (T^1_{c r_2 - c r_1 })^{\bigotimes (m-1)} S_1 
      (T^{1}_{c r_1})^{\bigotimes (m)} f^{\bigotimes (m)}, \mu_0^{\bigotimes (m-k)} \>  \nonumber \\ 
  \le && K \left[\prod_{i=1}^k \sup_{x \in \R^d}T^1_{c r_i}f(x)\right] \<(T^1_{c t})^{\bigotimes (m-k)}f^{\bigotimes (m-k)}, \mu_0^{\bigotimes (m-k)}\> 
  \nonumber \\ 
  = && K \left[\prod_{i=1}^k \sup_{x \in \R^d}T^1_{c r_i}f(x)\right] \sup_{0 < t \leq T} \<T^1_{ct}f,\mu_0 \>^{(m-k)}. \nonumber 
   \eeqlb 
   By Hypothesis $\ref{hyp:basicassumpGauss}$ and $(\ref{LSU})$, we get $\sup_{0\le t\le T}|M_k(t)|<\infty$ in the case 
   $f(w)=q_\epsilon^{1}(w,x)\ge0$. 
   
In the second case we apply the same argument to $f(w)=|\partial_p q_\epsilon^{1}(w,x)|$ and for simplicity we set $c=1$. 
Equation $(\ref{LSU})$ immediately allows us to bound $T^1_{r_i}f$ by $r_i^{1/2}T^1_{r_i}\varphi$.  
Since operator $\partial_p$ commutes with operator $T_{t}^1$, 
without loss of generality we restrict ourselves to the exponential jump times 
$\tau_i \sim \frac{1}{2}(m-i+1)(m-i)e^{-\frac{1}{2}(m-i+1)(m-i) }$ of the dual process. Fixing $t^* = r_i$ for some arbitrary $i$, 
by the dominated convergence theorem, we have
\beqlb \lab{5ES5}
 \sup_{x \in \R^d}T^1_{t^*}f(x) & \leq & \sup_{x \in \R^d}  \int_0^{\infty} K
  \frac{e^{- \lambda r}}{\sqrt{r}} \<\varphi_r(w- \xi), \varphi_{t^*}(x- \xi)d \xi \> dr \nonumber \\
   &  \leq &  K  \sup_{x \in \R^d}  | T_{t^*}^1 \partial_p \tilde{Q}^{\lambda}(w-x)|   \nonumber  \\
    &  =  &  K \sup_{x \in \R^d}  |  \partial_p T_{t^*}^1 \tilde{Q}^{\lambda}(w-x)|   \nonumber  \\
     & = &  K   \sup_{x \in \R^d}|  \partial_p  \tilde{Q}_{t^*}^{\lambda}(w-x)|   \nonumber  \\
     & \leq &  K   \sup_{z \in \R^d}|  \partial_p  \tilde{Q}_{t^*}^{\lambda}(z)|  <   \infty,
\eeqlb
where $\tilde{Q}_{t^*}^{\lambda}(z)$ is defined by
$\tilde{Q}_t^{\lambda}(z) := \int_0^{\infty} e^{- \lambda r} \varphi_{t+r}(z) dr$. 
Next
 \beqlb \lab{5ES6}
  & &  \sup_{0 < t \leq T} \sup_{w \in \R^d} \<T^1_{t}f, \mu_0\>  \nonumber \\
  & &  \leq  \sup_{0 < t \leq T} \sup_{w \in \R^d} \int_{\R^d} \int_{\R^d} \varphi_t(x-y) f(y,w) dy \mu_0(dx)  \nonumber   \\
  & & \mbox{By Fubini's Theorem}  \nonumber    \\
  & & \leq \sup_{0 < t \leq T} \sup_{w \in \R^d} \int_{\R^d} \int_{\R^d} \varphi_t(x-y) \mu_0(dx)f(y,w)dy   \nonumber   \\
  & &  \leq \sup_{0 < t \leq T} \sup_{y \in \R^d}\<\varphi_t(x,y), \mu_0(dx)\> \sup_{w \in \R^d} \int_{\R^d}f(y,w)dy.
\eeqlb

By Hypothesis $\ref{hyp:basicassumpGauss}$, the only term remaining to bound is 
 \beqlb \lab{5ES8}
  & &   \sup_{w \in \R^d} \<f(y,w), \lambda_0\> \leq \sup_{w \in \R^d} \< \int_{0}^{\infty}e^{- \lambda s} \frac{1}{s^{(d+1)/2}}\exp{\{
   - \frac{a_0 |w-y|^2}{s} \}} ds, \lambda_0(dy)\> \nonumber \\
  & & \leq K \sup_{w \in \R^d}\< \int_{0}^{\infty}e^{- \lambda s} \frac{1}{s^{1/2}} \varphi_s(w-y) ds, \lambda_0(dy)\> \nonumber \\
  & &  \leq K \int_{0}^{\infty}e^{- \lambda s} \frac{1}{s^{1/2}} ds < \infty.
  \eeqlb
Hence $\sup_{0\le t\le T}|M_k(t)|<\infty$ ensues. 
 \qed 
 
  \proof [Proof of Theorem \ref{ClaimAk}] 
 We explicit the proof of both (\ref{AME}) and (\ref{BME}) only for the base case $k=1$, 
 the general case ensuing from (\ref{ME}) in Corollary \ref{MCOM} in a similar fashion. 
More precisely, for any $T > 0$, $0 \leq t \leq T$, $\epsilon > 0$ and $\epsilon^{'} > 0$, we have
 \beqnn
\sup_{x \in \R^d} \sup_{0 \leq t \leq T} \E_{\mu_{0}} |\Lambda^{x, \epsilon}_t| 
&& \leq \sup_{x \in \R^d} \sup_{0 \leq s \leq T} \E_{\mu_0} \int_0^t \<q_{\epsilon}(x- \cdot), \mu_s\>ds \nonumber \\
&& \leq \sup_{x \in \R^d} \sup_{0 \leq s \leq T}  \int_0^t \<P_sq_{\epsilon}(x- \cdot), \mu_0\>ds  \nonumber \\
&& \leq \sup_{x \in \R^d} \sup_{0 \leq s \leq T} \<q_{s + \epsilon}(x- \cdot), \mu_0\>   < \infty 
 \eeqnn
while, using also the sharp estimation (\ref{formulaI}) for the proof of (\ref{BME}),
 \beqnn
&& \lim_{\epsilon \downarrow 0}\lim_{\epsilon^{'} \downarrow 0}\sup_{x \in \R^d}\sup_{0 \leq t \leq T} \E_{\mu_0}| \Lambda^{x, \epsilon}_t - \Lambda^{x, \epsilon^{'}}_t| \nonumber \\
&& \leq \lim_{\epsilon \downarrow 0} \lim_{\epsilon^{'} \downarrow 0} \sup_{x \in \R^d}\sup_{0 \leq t \leq T} \E_{\mu_0}| \int_0^t \<(q_{\epsilon}(x - \cdot) - q_{\epsilon^{'}}(x - \cdot)), \mu_s\> ds|  \nonumber \\
&& \leq \lim_{\epsilon \downarrow 0} \lim_{\epsilon^{'} \downarrow 0} \sup_{x \in \R^d}\sup_{0 \leq t \leq T} \E_{\mu_0}| \int_0^t \<c|\epsilon - \epsilon^{'}|^{\alpha/2}(k_1 q_{\epsilon}(x - \cdot) + k_2 q_{\epsilon^{'}}(x - \cdot)), \mu_s\> ds| = 0 \nonumber \\
 \eeqnn
holds for any $\alpha \in (0 , 1)$. 
\qed


  \section{Tanaka formula and local time of SDSM}\label{sec:Tanaka}
  \setcounter{equation}{0}

 The stage is now set for the proof of our main result (Theorem \ref{lt_th1}). We need three lemmata.
 The first one concerns the Laplace transform $Q^{\lambda}$ defined in (\ref{eqn:green}).
 \ble \lab{le_b}
  Assume Hypothesis $\ref{hyp:basicassumpElliptic}$ is satisfied. For any $\lambda > 0$ there holds:\\
   {\em (i)} For all $d\ge1$, we have $Q^{\lambda} \in L^1(\R^d)$ and
  $\partial_{x_i}Q^{\lambda} \in L^1(\R^d)$ for any $i \in \{1,2,\cdots,d\}$. \\
  {\em (ii)} For $d=1$, we also have $\partial_{x} Q^{\lambda} \in L^2(\R)$. \\
   {\em (iii)} For $d=1$, $2$ or $3$, we finally have $Q^{\lambda} \in L^2(\R^d)$. \\
 If Hypotheses $\ref{hyp:basicassumpGauss}$ and $\ref{hyp:basicassumpUniformInteg}$ 
 are also satisfied for $\mu_0 \in M_a(\R^d)$ with $a\ge0$, then
   there holds : \\
  {\em (iv)} For all $d\ge1$, we have $Q^{\lambda} \in L^1(\mu_0)$ and
  $\partial_{x_i}Q^{\lambda} \in L^1(\mu_0)$ for any $i \in \{1,2,\cdots,d\}$. \\
   {\em (v)} For $d=1$, $2$ or $3$, we finally have $Q^{\lambda} \in L^2(\mu_0)$.
 \ele

 The proof is technical and found in Subsection \ref{app:pf_le_b}. \\

Further, for each $\lambda > 0$ and $x \in \R^d$, $u(\cdot)=Q^{\lambda}(x-\cdot)$ solves equation
 $\left( - G_1  + \lambda \right) u = \delta_{x}$
 in the distributional sense, so the Green operator $Q^{\lambda}*\phi(x)=\int_{\R^d}\phi(y)Q^{\lambda}(x-y)dy$
 for Markov semigroup $P^1_t$, is a well-defined convolution for any $\phi\in C_b(\R^d)$ and solves
 \beqlb \lab{GPDE}
      (- G_1 + \lambda)u = \phi.
 \eeqlb

  \ble \lab{le_new1}
  Assume Hypotheses $\ref{hyp:basicassumpFilter}$, 
  $\ref{hyp:basicassumpElliptic}$, $\ref{hyp:basicassumpGauss}$ and $\ref{hyp:basicassumpUniformInteg}$ are satisfied.
  For either $d=1$, $2$ or $3$, the random field
 \beqlb \lab{DI0}
  \Xi_t(x):= \int_0^{t} \int_{\R^d} \<h_p(y-\cdot) \partial_pQ^{\lambda}(x-\cdot),
 \mu_s\> W(dy,ds)
 \eeqlb
is a square-integrable ${\cal F}_t$-martingale,
for every $\lambda > 0$ and $p\in\{1,2,\ldots,d\}$, with quadratic variation given by
\beqnn
  \<\Xi(x)\>_t
  = \int_0^{t} ds \int_{\R^d} \<h_p(y-\cdot)\partial_pQ^{\lambda}(x-\cdot), \mu_s\>^2 dy
 \eeqnn
 and satisfying $\sup_{x\in\R^d} \E_{\mu_0}  \<\Xi(x)\>_t  < \infty $ for every $t > 0$.
  \ele

  \proof
 Recalling the single particle transition density $q_t^{1}(0,\cdot)$ in the case where the starting position
 is the origin $0\in\R^d$, define the following perturbation of $Q^{\lambda}$ from (\ref{eqn:green}),
for every $\lambda > 0$ and $\epsilon > 0$:
\beqlb \lab{perturbation}
Q^{\lambda}_{\epsilon}(\cdot) := \int_{0}^{\infty}e^{- \lambda u} q_{u + \epsilon}^{1}(0,\cdot)du = e^{\lambda \epsilon}\int_{\epsilon}^{\infty}e^{- \lambda t} q_t^{1}(0,\cdot)dt
\eeqlb
and observe that $Q^{\lambda}_{\epsilon}\in C_b^{\infty}(\R^d)$  and
\beqlb \lab{basic1}
  \lim_{\epsilon \downarrow 0}|Q^{\lambda}_{\epsilon}(x) -  Q^{\lambda}(x)| \leq  \lim_{\epsilon \downarrow 0} \int_0^{\infty} e^{- \lambda u}|P^1_{\epsilon}*q_u^{1}(x) -
 q_u^{1}(x)|du = 0,
\eeqlb
for every $x \in \R^d \setminus \{y= (y_1, y_2, \cdots,y_d): y_1 y_2 \cdots y_d = 0\}$ 
as $\epsilon \downarrow 0 $. By (iv) and (v) of Lemma \ref{le_b} and the
dominated convergence theorem, there holds, for every $a\ge0$ and $\mu_0 \in M_a(\R^d)$,
\[
Q^{\lambda}_{\epsilon} \rightarrow Q^{\lambda} \hspace{1cm} \mbox{    in $L^{k}(\mu_0); k=1,2; d=1,2,3;$}
\]
and
\[
\partial_{x_i}Q^{\lambda}_{\epsilon} \rightarrow \partial_{x_i}Q^{\lambda} 
\hspace{1cm} \mbox{in  $L^1(\mu_0); d \geq 1; i=1, \cdots, d.$ }
\]

In fact, by (\ref{LSU}) we know that  
${\partial_{j}^{r}} {\partial_{k}^{s}} Q^{\lambda}_{\epsilon}\in L^p(\mu_0)\cap C_0(\R^d)$
also holds for any choice of $p\ge1$, $a\ge0$, $\mu_0\in M_a(\R^d)$,
$1\le j, k \le d$ and nonnegative integers $r$ and $s$ such that $0\le r, s\le2$.
By Corollary \ref{MCOM}, we also know that
${\partial_{j}^{r}} {\partial_{k}^{t}} Q^{\lambda}_{\epsilon}\in L^p(\mu_t)$
holds $\P_{\mu_0}$-almost surely for all $t\ge0$. Hence, by Hypothesis $\ref{hyp:basicassumpElliptic}$,
the integrand in (\ref{DI0}) is $\P_{\mu_0}$-almost surely finite for all $t\ge0$.
We need a bit more, namely integrability with respect to $W$, which we prove next.

With $\Xi_t^{\epsilon}(x)$ defined like $\Xi_t(x)$ in (\ref{DI0}) with $Q^{\lambda}$ replaced by $Q^{\lambda}_{\epsilon}$,
we begin by showing that $\{\Xi_t^{\epsilon}(x):t\ge0\}$ is itself a square-integrable ${\cal F}_t$-martingale for every $x\in\R^d$.
For any choice of $x,y\in\R^d$ let us denote by
$\Psi_{xy}^\epsilon\in C_b({(\R^d)^2})\cap L^1(\R^{2d})\cap L^1(\mu_0^2)$ the function
\[
\Psi_{xy}^\epsilon(w_1,w_2)
:= h_p(y-w_1) \partial_pQ^{\lambda}_{\epsilon}(x-w_1)\times h_p(y-w_2) \partial_pQ^{\lambda}_{\epsilon}(x-w_2).
\]
By Corollary \ref{MCOM}, we can apply the duality identity (\ref{dual})
to the second moment, to get
\beqlb \lab{DI1}
   \E_{\mu_0} \left[ \{\Xi_t^{\epsilon}(x)\}^2\right]
   & = & \E_{\mu_0} \int_0^t \int_{\r^d} \<\Psi_{xy}^\epsilon, \mu_s^2\> dy ds  \\
   & = & \int_0^t ds\int_{\r^d} dy\bigg[\<P^2_s\Psi_{xy}^\epsilon, \mu_0^2\> +
   \int_0^s \<P^1_{s-r}\Phi_{12}^2P^2_r\Psi_{xy}^\epsilon, \mu_0\> \gamma \sigma^2 dr\bigg], \nonumber
  \eeqlb
 where $P^1_s$ and $P^2_s$ are the one and two particle transition semigroups
 from (\ref{eqn:Semigroup}), with respective generators $G_1$ and $G_2$,
 and $\Phi_{12}^2$ is the diagonalisation mapping from (\ref{restriction}).

 Rewrite (\ref{LSU}) for every $(t,x,y)\in(0,T) \times (\R^d)^m \times (\R^d)^m$ with $T > 0$ as
 \beqlb \lab{LSU2}
   \left|\frac{\partial^{r}}{\partial t}\frac{\partial^{s}}{\partial y_{p}} q_t^{m}(x,y)\right|
   \leq \frac{a_1}{t^{(2r+s)/2}} \left(\frac{\pi}{a_2}\right)^{dm/2} \prod_{i=1}^m g_t(y_i-x_i),
  \eeqlb
 where $g_t$ is the centered gaussian density on $\R^d$ with $d$ independent coordinates
 and diagonal variance matrix with all entries equal to $t/2a_2$. 
 Use Hypothesis $\ref{hyp:basicassumpElliptic}$ to first take care of the term
 \beqlb \lab{DI10}
   & & \int_{\r^d} | h_p(y-w_1) | | h_p(y-w_2) | dy \leq \|h\|_\infty \|h\|_1 < \infty.
   \eeqlb 
   
 We can now bound $\Psi_{xy}^\epsilon$ by
 using the values $m=1$, $r=0$ and $s=1$ in (\ref{LSU2}), to get, for some constant $c=c(T,d,h) > 0$ and all $\epsilon > 0$,
 \beqlb \lab{Psibound}
 \int_{\r^d} | \Psi_{xy}^\epsilon(w_1,w_2) | dy
 \leq c \int_0^{\infty}\int_0^{\infty} e^{- \lambda (u+v)} \frac{1}{\sqrt{uv}} g_u(x-w_1)g_v(x-w_2) dudv.
  \eeqlb
  
 Therefore, using Fubini's theorem plus the values $m=2$ and $r=s=0$ in (\ref{LSU2}) this time, we get,
 with a new constant $c=c(T,d,h) > 0$,
 \beqlb \lab{DI8}
   \int_{\r^d} | P^2_r \Psi_{xy}^\epsilon(z_1,z_2) | dy
  & \leq & c \int_0^{\infty}\int_0^{\infty} e^{- \lambda (u+v)} \frac{1}{\sqrt{uv}} dudv  \nonumber  \\
  & & \times
  \int_{\r^d}\int_{\r^d} g_u(x-w_1)g_v(x-w_2)g_r(w_1-z_1)g_r(w_2-z_2) dw_1dw_2  \nonumber  \\
   & \leq & c \int_0^{\infty}\int_0^{\infty} e^{- \lambda (u+v)} \frac{1}{\sqrt{uv}} g_{r+u}(x-z_1)g_{r+v}(x-z_2) dudv.
  \eeqlb

By the same reasoning one finally obtains, with a further value for $c=c(T,d,h) > 0$,
 \beqlb \lab{DI9}
 \int_{\r^d} | P^1_{s-r}\Phi_{12}^2P^2_r \Psi_{xy}^\epsilon(\zeta) | dy
 & \leq & c \int_0^{\infty}\int_0^{\infty} e^{- \lambda (u+v)} \frac{1}{\sqrt{uv}} dudv  \nonumber  \\
 & & \times
 \int_{\r^d} g_{r+u}(x-z)g_{r+v}(x-z) q_{s-r}^{1}(\zeta,z) dz.
  \eeqlb

Hypothesis $\ref{hyp:basicassumpGauss}$ now yields a new $c=c(T,d,h,a,\mu_0)$ such that
  \beqlb \lab{DI12}
  & &  \int_0^t ds \int_{\r^d} dy \int_0^s  | \< P^1_{s-r}\Phi_{12}^2P^2_r\Psi_{xy}^\epsilon, \mu_0 \>  | dr
  = \int_0^t dr \int_{\r^d} dy \int_r^t  | \< P^1_{s-r}\Phi_{12}^2P^2_r\Psi_{xy}^\epsilon, \mu_0 \>  | ds \nonumber  \\
  & & \leq c \int_0^t dr \int_0^{\infty}\int_0^{\infty} e^{- \lambda (u+v)} \frac{1}{\sqrt{uv}} dudv
  \left[ \int_{\r^d} g_{r+u}(x-z)g_{r+v}(x-z)dz \right]  \nonumber  \\
  & & \leq ca_1\chi_{d} (t),
  \eeqlb
which is finite when $d\le3$, since (see Subsection \ref{app:pf_DI13} for the proof of this upper bound)
 \beqlb \lab{DI13}
 \chi_{d}(t) \equiv  \int_0^{\infty} \int_0^{\infty} e^{- \lambda (u+v)} \frac{1}{\sqrt{uv}}
  \int_0^t \frac{1}{(\sqrt{2r +u + v})^{d}} dr du dv
 \le  (t+1)\pi \sqrt{\frac{\pi}{\lambda}}.
 \eeqlb
This takes care of the second term in (\ref{DI1}); the first one is done similarly by way of (\ref{DI8}).

This proves that $\{\Xi_t^\epsilon(x):t\ge0\}$ is a square-integrable martingale for every $x\in\R^d$
and $\epsilon > 0$, as long as $d\le3$; moreover, for every $T > 0$ and $\lambda > 0$, there holds
\beqlb \lab{SquareInteg}
   \sup_{\epsilon>0}\sup_{x\in\R^d}\sup_{0 \leq t \leq T} \E_{\mu_0} \left[ \Xi_t^{\epsilon}(x)\right]^2 < \infty.
\eeqlb 

Using similar ideas we get
\beqlb \lab{DI1a}
   \sup_{x\in\R^d}\sup_{0 \leq t \leq T}
   \E_{\mu_0} \left[ \left| \Xi_t^{\epsilon}(x)-\Xi_t^{\delta}(x)\right|^2 \right] 
    \leq c \epsilon, 
  \eeqlb
for some new constant $c=c(T,d,h,a,\mu_0,a_1,a_2,\gamma,\sigma, \lambda) > 0$ 
 (this one dependent on $\lambda$ but independent of $\epsilon$) and 
 all $0 < \delta < \epsilon < 1$. Details are found in Subsection \ref{app:pf_DI1a}. 
Picking any sequence $\epsilon_n$ decreasing to $0$, we get that
$\{\Xi_t^0(x):=\lim_{\epsilon_n\rightarrow0}\Xi_t^{\epsilon_n}(x):t\ge0\}$ exists $\P_{\mu_0}$-almost surely and 
is a square-integrable martingale for every $x\in\R^d$. 
We define $\Xi_t$ in (\ref{DI0}) as this unique limit $\Xi_t:=\Xi_t^0$, 
in both of the $\P_{\mu_0}$-almost sure and $L^2(\P_{\mu_0})$ sense. 
The convergence in $L^2(\P_{\mu_0})$ of the sequence 
$\{\Xi_t^{\epsilon_n}(x)\}$, for each fixed $t\ge0$ and $x\in\R^d$,
together with Corollary \ref{MCOM} plus Hypotheses $\ref{hyp:basicassumpFilter}$, 
 $\ref{hyp:basicassumpElliptic}$ and $\ref{hyp:basicassumpGauss}$,
imply that (\ref{DI1}) holds also in the limiting case $\epsilon=0$. 

Since the right hand side of (\ref{DI8}), (\ref{DI12}), and (\ref{DI13}) do not dependent on $\epsilon$,  
by the dominated convergence theorem, we also have
$\sup_{x\in\R^d}\sup_{0 \leq t \leq T} \E_{\mu_0} \left[ \Xi_t(x)\right]^2  < \infty$.
\qed

 \ble \lab{SFT}
   Assume Hypotheses $\ref{hyp:basicassumpFilter}$, 
  $\ref{hyp:basicassumpElliptic}$, $\ref{hyp:basicassumpGauss}$ and $\ref{hyp:basicassumpUniformInteg}$ are satisfied. 
  For any $t>0$ and $\phi \in C^{\infty}_c(\R^d)$, we have $\P_{\mu_0}$-almost surely
  \beqlb \lab{ex1}
  & & \int_{\R^d} \phi(x) \int_0^{t} \int_{\R^d} \sum_{p=1}^{d} \<h_p(y- \cdot)
    \partial_p Q^{\lambda}(x - \cdot), \mu_s\>W(dy,ds) dx \nonumber \\
  & &  =  \int_0^{t} \int_{\R^d} \int_{\R^d} \phi(x) \sum_{p=1}^{d} \<h_p(y- \cdot)
    \partial_p Q^{\lambda}(x - \cdot), \mu_s\>dx W(dy,ds)
 \eeqlb
 and
 \beqlb \lab{ex2}
   & & \int_{\R^d} \phi(x) \int_0^{t} \int_{\R^d}Q^{\lambda}(x-y)M(dy,ds)dx  \nonumber \\
   & & =  \int_0^{t} \int_{\R^d} \int_{\R^d} \phi(x) Q^{\lambda}(x-y)dx M(dy,ds),
   \eeqlb
 where $W$ is a space-time white noise and $M$ is an orthogonal worthy martingale measure with
 \[
  \< \int_{\R^d} \phi(x-y) M(dyds)\> = \gamma \sigma^2 \int_{\R^d} \phi(x)^2 dx ds.
 \]
 \ele
 \proof 
  Both equalities follow from a stochastic version of Fubini's theorem due to Walsh \cite{Walsh86} (see his Theorem 2.6 on p. 296).
Since the covariance measure of $M(dx,ds)$ is deterministic, (\ref{ex1}) and (\ref{ex2}) are handled similarly; 
we only prove the more difficult (\ref{ex1}).  By Corollary 2.8 of Walsh \cite{Walsh86}, we only need to check 
 \beqlb \lab{cond1}
  & & \hspace{-1cm} \E_{\mu_0} \int_{\R^d} \int_{\R^d} \int_0^t  \<\sum_{p=1}^{d}|h_p(\xi - y)
    \partial_p Q^{\lambda}(x - y)|, \mu_s(dy)\>  \nonumber  \\
      & & \cdot \< \sum_{q=1}^{d} | h_q(\xi - y)
    \partial_q Q^{\lambda}(x - y)|, \mu_s(dy)\> d\xi ds \nu(dx) < \infty,
 \eeqlb
 for any finite measure on ${\cal B}(\R^d)$ of the form $\nu(dx) = \phi(x) dx$ with positive $\phi \in C_c(\R^d)$. 
 Once again write $g_r$ for the centered gaussian density on $\R^d$ with $d$ independent coordinates
 and diagonal variance matrix with all entries equal to $r/2a_2$. Let
 \[
  f(x-y) := \int_0^{\infty}
  \frac{e^{- \lambda r}}{\sqrt{r}} g_r(x - y)dr
 \] 
and observe
  \beqlb \lab{cond2}
  & & \hspace{-1cm} \E_{\mu_0} \int_{\R^d} \int_{\R^d} \int_0^t  \<\sum_{p=1}^{d}|h_p(\xi - y)
    \partial_p Q^{\lambda}(x - y)|, \mu_u(dy)\>  \nonumber  \\
      & & \cdot \< \sum_{q=1}^{d} | h_q(\xi - y)
    \partial_q Q^{\lambda}(x - y)|, \mu_u(dy)\> d\xi du \nu(dx) \nonumber  \\
  & &  \leq k_1 \sum_{p,q=1}^{d} \rho_{pq}(0)  \E_{\mu_0} \int_{\R^d}  \int_0^t  \<\int_0^{\infty}
  \frac{e^{- \lambda r}}{\sqrt{r}} g_r(x - y)dr, \mu_u(dy)\>  \nonumber  \\
      & & \cdot \<\int_0^{\infty}
  \frac{e^{- \lambda s}}{\sqrt{s}}g_s(x - y)ds, \mu_u(dy)\> du \nu(dx)  \nonumber  \\
  & &  \leq k_1 \sum_{p,q=1}^{d} \rho_{pq}(0) \int_0^t \int_{\R^d}  \E_{\mu_0}  \<f(x-y), \mu_u(dy)\>^2 du \nu(dx).
 \eeqlb
 We already proved the finiteness of
 \[
 \sup_{0 < u \leq T}\sup_{x \in \R^d} \E_{\mu_0}\<f(x-y), \mu_u(dy)\>^2 < \infty.
 \]
 Thus, (\ref{cond1}) holds. This proves (\ref{ex1}).
 \qed  \\
 
 \proof [Proof of Theorem \ref{lt_th1}]
 First, since $u=Q^{\lambda}$ is the fundamental solution of
 \[
      (\lambda - G_1)u = 0,
 \]
 we have
 \[
   (\lambda - G_1 ) Q^{\lambda}_{\epsilon}(x)= (\lambda - G_1 ) Q^{\lambda}*q_{\epsilon}^{1}(x) = q_{\epsilon}^{1}(x).
 \]
 
 For every $ \phi \in C_c(\R^d)$ and any given $t > 0$, by definition of $\Lambda^x_t $ we have
 \beqlb \lab{App0}
 && \lim_{\epsilon \downarrow 0}\E_{\mu_0}|\int_{\R^d} \phi(x) \Lambda^x_tdx - \int_{\R^d}\phi(x) \Lambda^{x, \epsilon}_t dx| \nonumber \\
 & & \leq \lim_{\epsilon \downarrow 0 }\sup_{x \in \R^d} \sup_{0 \leq t \leq T} \E_{\mu_0}|\Lambda^x_t - \Lambda^{x, \epsilon}_t|\int_{\R^d}\phi(x) dx  = 0.
 \eeqlb
Then, for every $0 \leq t \leq T$, $x \in \R^d$ and $\omega$ outside of the $\P_{\mu_0}$-null set $N$ specified by the statement of 
 Lemma \ref{SFT}, a stochastic version of Fubini's theorem, there holds
 \beqlb \lab{App1}
 & & \int_{\R^d} \phi(x) \Lambda^{x}_t dx = \lim_{\epsilon \downarrow 0}\int_{\R^d} \phi(x) \Lambda^{x, \epsilon}_t dx \nonumber \\
 & & = \lim_{\epsilon \downarrow 0} \int_{0}^{t}\<\int_{\R^d}(-G_1 + \lambda) Q^{\lambda}_{\epsilon}(x- \cdot) \phi(x)dx, \mu_s\> ds \nonumber \\
 &&  =   \lim_{\epsilon \downarrow 0} \int_{0}^{t}\<\int_{\R^d}q_{\epsilon}^{1}(x- \cdot) \phi(x)dx, \mu_s\> ds          \nonumber \\
  && = \lim_{\epsilon \downarrow 0} \int_{0}^{t}\<P_{\epsilon}*\phi(\cdot), \mu_s\> ds          =\int_0^t\<\phi, \mu_s\>ds.
\eeqlb
 This proves that $\Lambda^x_t$ is the SDSM local time.
 Since $Q^{\lambda} \ast \phi(\cdot) \in C_b^{\infty}(\R^d)$, by Lemma \ref{SFT}, we have 
 \beqlb \lab{App1}
 & & \int_{\R^d} \phi(x) \Lambda^{x}_t dx = \int_{0}^{t}\<(-G_1 + \lambda) \int_{\R^d}Q^{\lambda}(x- \cdot) \phi(x)dx, \mu_s\> ds  \nonumber \\
  & & = \int_{0}^{t}\<(-G_1 + \lambda) Q^{\lambda}* \phi(\cdot), \mu_s\> ds \nonumber \\
  & &  =  \<Q^{\lambda}* \phi(\cdot), \mu_0\> - \<Q^{\lambda}* \phi(\cdot),
    \mu_{t}\>
    + \lambda \int_0^{t} \<Q^{\lambda}* \phi(\cdot), \mu_s\>ds \nonumber \\
    & &  + \sum_{p=1}^{d}\int_0^{t} \int_{\R^d}\<h_p(y- \cdot)
    \partial_p Q^{\lambda}* \phi(\cdot), \mu_s\>W(dy,ds) \nonumber \\
   & &  + \int_0^{t} \int_{\R^d}Q^{\lambda}* \phi(y)M(dy,ds) \nonumber \\
   & &  = \int_{\R^d} \phi(x) \bigg\{\<Q^{\lambda}(x - \cdot), \mu_0\> - \<Q^{\lambda}(x- \cdot),
    \mu_{t}\>
    + \lambda \int_0^{t} \<Q^{\lambda}(x - \cdot), \mu_s\>ds \nonumber \\
    & &  + \sum_{p=1}^{d}\int_0^{t} \int_{\R^d}\<h_p(y- \cdot)
    \partial_p Q^{\lambda}(x - \cdot), \mu_s\>W(dy,ds) \nonumber \\
   & &  + \int_0^{t} \int_{\R^d}Q^{\lambda}(x-y)M(dy,ds)\bigg\}dx
   \eeqlb
 $\P_{\mu_0}$-almost surely, jointly for all $\phi \in C_c(\R^d)$. So (\ref{TanakaII}) holds and we are done.
 \qed


\section{Proofs of technical results}\label{app:ProofsOfLemmas}
\setcounter{equation}{0}

\subsection{Proof of uniform ellipticity of operator ${G}_{m}$ in (\ref{eqn:Gn})}\label{app:pf_ellip1}
 \proof 
To check the uniform ellipticity of $(\Gamma^{ij}_{pq})$, let $\xi_i= (\xi_{i1}, \cdots, \xi_{id})^{T}$
 denote an arbitrary column vector in $\R^d$ and
 $\Gamma := (\Gamma^{ij}_{pq}(x_1, \cdots, x_m))_{1 \leq i,j \leq m, 1 \leq p,q \leq d}$. Since
 \beqlb \lab{ellip1}
 & & (\xi_1^{T}, \cdots, \xi_m^{T} ) \Gamma  \left( \ba{c}
   \xi_1  \\
    \vdots \\
   \xi_{m} \ea \right) = \sum_{i,j=1}^{m} \sum_{p,q=1}^{d} \xi_{ip}
   \Gamma^{ij}_{pq}(x_1, \cdots, x_m) \xi_{jq} \nonumber \\
 &&  = \sum_{i=1}^{m} \sum_{p,q=1}^{d}\left[\xi_{ip}(a_{pq}(x_i) +
 \rho_{pq}(x_i,x_i)) \xi_{iq} \right]
  + \sum_{i,j=1, i \neq j}^{m} \sum_{p,q=1}^{d}\xi_{ip} \rho_{pq}(x_i,x_j)
 \xi_{jq} \nonumber \\
 & & = \sum_{i=1}^{m} \bigg[\sum_{r=1}^{d}
 \big(\sum_{p=1}^d \xi_{ip} c_{pr}(x_i) \big)^2
 + \int_{\R^d}\big(\sum_{p=1}^d \xi_{ip} h_{p}(u-x_i)\big)^2du  \bigg]
 \nonumber \\
 & & + \sum_{i,j=1, i \neq j}^{m} \int_{\R^d} \big( \sum_{p=1}^d
 \xi_{ip} h_{p}(u-x_i) \big)\big( \sum_{p=1}^d
 \xi_{jq} h_{q}(u-x_j) \big) du \nonumber \\
 & & = \sum_{i=1}^{m} \sum_{r=1}^{d}
 \big(\sum_{p=1}^d \xi_{ip} c_{pr}(x_i) \big)^2
 + \int_{\R^d}\bigg[ \sum_{i=1}^{m} \big(\sum_{p=1}^d \xi_{ip}
 h_{p}(u-x_i)\big)\bigg]^2du \geq 0 , 
 \eeqlb
 by the uniform ellipticity assumption of $(a_{pq})_{1 \leq p,q \leq d}$
  there exists a positive real number $\epsilon > 0$ such that
 for each $1 \leq i \leq m$
 \beqlb \lab{u_e}
 \sum_{r=1}^{d}
 \big(\sum_{p=1}^d \xi_{ip} c_{pr}(x_i) \big)^2=\sum_{p,q=1}^{d}\left[\xi_{ip}a_{pq}(x_i) \xi_{iq} \right]
 \geq \epsilon | \xi_i|^2,
 \eeqlb
 where $| \xi_i|= \sqrt{\xi_{i1}^2 + \cdots + \xi_{id}^2}$ . The
 uniform ellipticity of $\Gamma$ follows. In the left hand side of the last inequality of (\ref{ellip1}), the first term is associated with the coefficients of the individual motions and the second term is associated with the coefficients of the common random environment. 
 \qed

\subsection{Proof of uniformity of Aronson bounds in (\ref{Aronsonbounds})}\label{app:pf_jointparametrization}
 \proof 
 Since mapping $c\mapsto c^{d/2}\varphi_{cs}(y-x)$ is increasing for any $x, y \in \R^d$ and $s > 0$, 
 all upper and lower Aronson bounds are respectively bounded 
 uniformly in $d\in\{1, 2, \ldots, n\}$ above and below, via  
$$  
\bigg(\frac{c_*}{c^*}\bigg)^{n/2}  \alpha_* \cdot \varphi_{c_*s}(y-x) \leq
A^*(d) \cdot \varphi_{c(d)s}(y-x) \leq \bigg(\frac{c^*}{c_*}\bigg)^{n/2}  \alpha^* \cdot \varphi_{c^*s}(y-x)
$$
  where $c=c(d)$ expresses the dependency on the dimension, $c^*=\max\{c(1), c(2), \ldots, c(n)\}$ and 
  $c_*=\min\{c(1), c(2), \ldots, c(n)\}$, while $\alpha^*=\max\{A^*(1), A^*(2), \ldots, A^*(n)\}$ and finallly 
  $\alpha_*=\min\{A^*(1), A^*(2), \ldots, A^*(n)\}$. 
   \qed
   
\subsection{Proof of rescaling properties in Section \ref{sec:dualConst}}\label{app:pf_iam}

\proof [Proof of statement (\ref{iam})]
 Clearly it suffices to prove it in the case $m=1$. For every $d\ge1$ we see at once that $I_a\in C_b(\R^d)$ holds.
 The symmetry of $|x|^2$ in its coordinates also
 allows for a reordering of the partial derivatives.
 Using repeatedly the formula $\partial_{x_i} I_a(x)=-ax_i I_{a+2}(x)$, by induction on $k\ge0$ we have
 \begin{equation*}\partial_{x_1}^k I_a(x)=\sum_{\gamma=0}^k H_{\gamma,k}(x_1)I_{a+k+\gamma}(x)
\end{equation*}
for every $d\ge1$, where $H_{\gamma,k}$ is a real polynomial in $x_1$ of order at most ${\gamma}$.
Therefore all these partial derivatives are finite sums of bounded continuous functions.
Note that $H_{\gamma,k}$ depends on parameter $a\ge0$ but not on dimension $d\ge1$.
This also finishes the proof in the base case $d=1$. The next stage yields
 \begin{equation*}\partial_{x_2}^{\ell}\partial_{x_1}^k I_a(x)
 =\sum_{{\beta}=0}^{\ell}\sum_{{\gamma}=0}^k H_{\gamma,k}(x_1)H_{\beta,\ell}(x_2)I_{a+\ell+k+\beta+\gamma}(x)
\end{equation*}
for every $k\ge1$, $\ell\ge1$ and $d\ge1$, again all bounded and continuous. Case $d=2$ is also finished.
Another iteration finishes the proof,
this time through a multiple induction on the integer vectors of the form $(d,k_1,\ldots,k_d)$
with $d\ge1$ and every $k_i\ge0$.
\qed

\proof [Proof of statement (\ref{iamgm})]
Using the above calculations, one can also show that
${\cal I}_{a, m}^{-1} G_m {\cal I}_{a, m} \in C_b((\R^d)^{m})$ holds for every $m\ge1$.
Indeed, since all coefficients $\Gamma_{pq}^{ij}$ of $G_m$ in (\ref{eqn:Gn}) are bounded and continuous
by Hypothesis \ref{hyp:basicassumpElliptic}, with $G_m$ involving only partial derivatives of
order $1$ or $2$, it suffices to prove this statement for $m\le2$.
In the base case $m=1$, $I_a^{-1} G_1 I_a \in C_b(\R^d)$ holds because both
$I_a^{-1}(x) \partial_{x_1}^2 I_a(x)$
and $I_a^{-1}(x) \partial_{x_2}\partial_{x_1} I_a(x)$ are bounded and continuous.
For $m=2$ the only new terms take the form
$I_a^{-1}(x)I_a^{-1}(y)\partial_{x_1} I_a(x)\partial_{y_1}I_a(y)$, a product of two terms
already covered in case $m=1$.
\qed

\proof [Proof of Lemma \ref{lea}]
Proposition 1.1.5b of Ethier and Kurtz \cite{EthierKurtz86} ensures that both $P_t^{m}f\in{\cal D}(G_m)$ and 
$ \frac{\partial}{\partial t} P_t^{m} f = G_m P_t^{m} f = P_t^{m} G_m f$ hold 
for every choice of $f\in{\cal D}(G_m)$ and $t\ge0$. 
The result is therefore true in the case $a=0$ under Hypothesis \ref{hyp:basicassumpElliptic}, 
so let us proceed with parameter values $a > 0$ in mind. All we need to show is that
both the following integrability conditions hold, for any $f \in {\cal D}_a(G_m)$: 
$ \sup_{0 \leq t \leq T} \| {\cal I}_{a, m}^{-1} P^{m}_t f \|_{\infty} < \infty$ and
$ \sup_{0 \leq t \leq T} \| {\cal I}_{a, m}^{-1} P_t^{m} G_m f \|_{\infty} < \infty$. 
It suffices to show their validity when $f={\cal I}_{a, m}$
because the linearity of $P_t^{m}$ implies both 
$ | {\cal I}_{a, m}^{-1} P^{m}_t f (x)| \le {\cal I}_{a, m}^{-1} P^{m}_t {\cal I}_{a, m}(x)
 \cdot \| {\cal I}_{a, m}^{-1} f \|_{\infty}$  and  
$ | {\cal I}_{a, m}^{-1} P^{m}_t G_m f (x)| \le {\cal I}_{a, m}^{-1} P^{m}_t {\cal I}_{a, m}(x)
 \cdot \| {\cal I}_{a, m}^{-1} G_m f \|_{\infty}$, for every $x\in(\R^d)^m$,  
since we already assume $\| {\cal I}_{a, m}^{-1} f \|_{\infty} < \infty$ and 
$\| {\cal I}_{a, m}^{-1} G_m f \|_{\infty} < \infty$ by choosing $f \in {\cal D}_a(G_m)$.

We first prove that 
$ \sup_{0 \leq t \leq T} \| {\cal I}_{a, m}^{-1} P^{m}_t {\cal I}_{a, m} \|_{\infty} < \infty$ holds.  
Since the gaussian density integrates to $1$ for every mean $x\in(\R^d)^m$, 
 as in $\< g_t^{m}(x,\cdot), \lambda_0 \>=1$ with  
  \beqlb \lab{Gauss}
     g_t^{m}(x,y)= \left(\frac{a_2}{\pi t} \right)^{dm/2}\exp{\left\{- a_2 \left( \frac{|y-x|^{2}}{t} \right)\right\}},
   \eeqlb
the application of inequality (\ref{LSU}) with $r=s=0$ yields upper bound
\beqnn
  {\cal I}_{a, m}^{-1}(x) P^{m}_t {\cal I}_{a, m} (x) 
   & = & {\cal I}_{a, m}^{-1}(x) \int_{(\R^d)^m} {\cal I}_{a, m} (y) q_t^{m}(x,y)dy \\
   & \le & a_1\left(\frac{\pi}{a_2} \right)^{dm/2} {\cal I}_{a, m}^{-1}(x) \int_{(\R^d)^m} {\cal I}_{a, m} (y) g_t^{m}(x,y)dy \\
   & = & a_1\left(\frac{\pi}{a_2} \right)^{dm/2} 
   \prod_{i=1}^m\left(I_a^{-1}(x_i)\int_{\R^d}I_a (y_i) g_t^{1}(x_i,y_i)dy_i\right)
 \eeqnn
for every $x\in(\R^d)^m$ and $t\in(0,T)$. 
Therefore this part of the proof will be complete once we show that 
there is a positive constant $C=C(a,d,T)$ such that (when $m=1$)
 \beqlb \lab{mainconstant}
\sup_{0 \leq t \leq T} \sup_{x\in\R^d}\left(I_a^{-1}(x)\int_{\R^d}I_a (y) g_t^{1}(x,y)dy\right)\le C
 \eeqlb
holds for all $a\ge0$. 
Any standard gaussian random variable $Z\sim N(0,1)$ satisfies the inequality 
$P(|Z|\ge z)\le e^{-z^2/2}$ for every $z\ge0$, so we have, 
for every $z\ge0$ and $t\ge0$,
\beqnn
 I_a^{-1}(x)\int_{\{y:|y-x|\ge z|x|\}}I_a (y) g_t^{1}(x,y)dy
 & \le & I_a^{-1}(x)\exp{\left\{- \frac{a_2 z^{2}|x|^{2}}{t}\right\}} \\
 & \le & \max\left\{2^{a/2},\left(\frac{at}{ea_2z^2}\right)^{a/2}\right\}, 
 \eeqnn
using $I_a \le1$ in the first inequality and then 
separating as to whether or not $|x|\le1$, to get the second one; 
restricting to $z\in(0,1/2)$ both ensures $|y|\ge (1-z)|x|$ 
on the complement of set $\{y:|y-x|\ge z|x|\}$ and makes 
$I_a^{-1}(x)I_a ((1-z)x)$ into an increasing function in $|x|$, so that
\[ 
I_a^{-1}(x)\int_{\{y:|y-x| < z|x|\}}I_a (y) g_t^{1}(x,y)dy
 \le I_a^{-1}(x)I_a ((1-z)x)\le(1-z)^{-a}.
\]
These bounds together yield the existence of such a $C$, with $C\ge2$ for all $a\ge0$. 

Finally, for every $x\in(\R^d)^m$ and $t\in(0,T)$, we have
\[ 
| {\cal I}_{a, m}^{-1} P_t^{m} G_m {\cal I}_{a, m} | (x) 
= | {\cal I}_{a, m}^{-1} P_t^{m} ({\cal I}_{a, m} {\cal I}_{a, m}^{-1} G_m {\cal I}_{a, m}) (x) | 
\le \| {\cal I}_{a, m}^{-1} G_m {\cal I}_{a, m} \|_{\infty} \cdot 
{\cal I}_{a, m}^{-1} P_t^{m} {\cal I}_{a, m} (x),
\]
so that  
$ \sup_{0 \leq t \leq T} \| {\cal I}_{a, m}^{-1} P_t^{m} G_m {\cal I}_{a, m} \|_{\infty} < \infty$ 
holds as well.
\qed

\proof [Proof that the mappings $\Phi_{ij}^m$ in (\ref{restriction}) are well-defined]
 We need to prove that $\Phi_{ij}^m$ maps ${\cal D}_a(G_m)$ into ${\cal D}_a(G_{m-1})$ for every $m\ge2$. 
 First notice the equalities 
 \beqnn
  & & {\cal I}_{a, m-1}^{-1}(y_1, \cdots, y_{j-1}, y_{j+1}, \cdots, y_m) 
  f(y_{1}, \cdots, y_{j-1},y_{i},y_{j+1},\cdots, y_{m})  \\
 = & & I_a(y_{i})I_a^{-1}(y_{i}){\cal I}_{a, m-1}^{-1}(y_1, \cdots, y_{j-1}, y_{j+1}, \cdots, y_m) 
  f(y_{1}, \cdots, y_{j-1},y_{i},y_{j+1},\cdots, y_{m})  \\
  = & & I_a(y_{j})I_a^{-1}(y_{j}){\cal I}_{a, m-1}^{-1}(y_1, \cdots, y_{j-1}, y_{j+1}, \cdots, y_m) 
  f(y_{1}, \cdots, y_{j-1},y_{j},y_{j+1},\cdots, y_{m}) \\
 =  & & I_a(y_{j}){\cal I}_{a, m}^{-1}(y_1, \cdots, y_{j-1}, y_{j}, y_{j+1}, \cdots, y_m) 
  f(y_{1}, \cdots, y_{j-1},y_{j},y_{j+1},\cdots, y_{m}), 
 \eeqnn
  the second one holding whenever variable $y_j$ is set to the value of $y_i$. Taking the supremum over 
  all $(y_1, \cdots, y_{j-1}, y_{j}, y_{j+1}, \cdots, y_m)\in(\R^d)^m$ in the top and bottom lines yields
 \beqlb \lab{shrink}
       || {\cal I}_{a, m-1}^{-1} \Phi_{ij}^mf ||_{\infty}\le || I_a ||_{\infty} \cdot || {\cal I}_{a, m}^{-1} f ||_{\infty}
       \le || {\cal I}_{a, m}^{-1} f ||_{\infty},  
 \eeqlb
so that $ || {\cal I}_{a, m}^{-1} f ||_{\infty} < \infty$ implies $ || {\cal I}_{a, m-1}^{-1} \Phi_{ij}^mf ||_{\infty} < \infty$, 
which is the first requirement. By the same reasoning, it follows that 
$ || {\cal I}_{a, m}^{-1} D_m f ||_{\infty} < \infty$ implies $ || {\cal I}_{a, m-1}^{-1} \Phi_{ij}^m D_mf ||_{\infty} < \infty$, 
for any differential operator $D_m$ of order at most $2$ (such as $G_m$, for example). 
The second requirement is $|| {\cal I}_{a, m-1}^{-1} G_{m-1} \Phi_{ij}^mf ||_{\infty} < \infty$, which follows at once upon
noting that the chain rule yields
$\partial/\partial_{x_{iq}}\Phi_{ij}^mf=\Phi_{ij}^m(\partial/\partial_{x_{iq}}f+\partial/\partial_{x_{jq}}f)$ 
while $\partial/\partial_{x_{jq}}\Phi_{ij}^mf=0$
and $\partial/\partial_{x_{kq}}\Phi_{ij}^mf=\Phi_{ij}^m(\partial/\partial_{x_{kq}}f)$ 
for every $k\not\in\{i,j\}$ and every $q$. Handling the second derivatives similarly one obtains
an operator $D_m$ such that $G_{m-1} \Phi_{ij}^mf = \Phi_{ij}^m D_mf$.
  \qed 

\proof [Proof of Lemma \ref{lea2}]
Given $J_0=m\ge1$ and $Y_0\in{\cal D}_a(G_m)$, rewrite the expression for process $Y$ in (\ref{Yprocess})
by replacing each $S_k$ with ${\cal I}_{a, k-1}{\cal I}_{a, k-1}^{-1}S_k{\cal I}_{a, k}{\cal I}_{a, k}^{-1}$ and reading
the new expression from right to left through the natural triplets thus formed, to get
\[
{\cal I}_{a, J_t}^{-1} Y_t = \left({\cal I}_{a, J_t}^{-1} P^{J_{\tau_k}}_{t-\tau_k} {\cal I}_{a, J_{\tau_{k-1}}}\right)
 \left({\cal I}_{a, J_{\tau_{k-1}}}^{-1} S_k {\cal I}_{a, J_{\tau_{k-2}}}\right) \cdots
 \left({\cal I}_{a, J_{\tau_{1}}}^{-1} S_1 {\cal I}_{a, J_{\tau_{0}}}\right)
 \left({\cal I}_{a, J_{0}}^{-1} P^{J_0}_{\tau_1}Y_0\right),
\]
for $\tau_k \le t < \tau_{k+1}$ and $0\le k\le J_0-1$.
By inequality (\ref{shrink}), all the triplets in $S$ are contractions and will not affect the overall uniform bound.
By the calculations leading to and including (\ref{mainconstant}), there exists a positive constant
$C=C(a,d,m,T)$ such that, for every $k\ge1$ and $f\in{\cal D}_a(G_k)$,
\[
  \sup_{0 \le t < \tau_{k+1}-\tau_k} || {\cal I}_{a, k}^{-1} P^{k}_t f ||_{\infty}
  \le \sup_{0 \le t < \tau_{k+1}-\tau_k} || {\cal I}_{a, k}^{-1} P^{k}_t  {\cal I}_{a, k} ||_{\infty}
   \cdot || {\cal I}_{a, k}^{-1} f ||_{\infty}
   \le C^k || {\cal I}_{a, k}^{-1} f ||_{\infty} < \infty.
 \]
This bound, when applied from right to left along any trajectory of rewritten process $Y$,
to the resulting triplets of the form ${\cal I}_{a, k}^{-1}P^{k}_t{\cal I}_{a, k}$,
yields the upper bound
\[
 \sup_{0 \leq t \leq T} || {\cal I}_{a, J_t}^{-1} Y_t ||_\infty
 \le C^{1+2+\ldots+m} || {\cal I}_{a, m}^{-1} Y_0 ||_\infty,
\]
for every $m\ge1$, since $C\ge2$ implies $C^{1+2+\ldots+\ell}\le C^{1+2+\ldots+m}$ for every $1\le \ell\le m$.
 \qed

\subsection{Proof of Lemma \ref{le_b}}\label{app:pf_le_b}
 \proof
 For any $d\ge1$ and $\lambda > 0$, $Q^{\lambda}(x)$
 is integrable as a direct consequence of (\ref{LSU}) with $m=1$, $y=0$, $r=0$ and $s=0$:
 there are constants $a_1 > 0$ and $a_2 > 0$ such that
  \beqnn
    \int_{\R^d} |Q^{\lambda}(x)| dx
    \leq  \int_0^{\infty}e^{-\lambda t} a_1
   \prod_{i=1}^{d} \left(\int_{\R} \frac{1}{\sqrt{t}}
   \exp\{ - a_2 \frac{x_i^2}{t}\}dx_i \right)dt
   = \left( \frac{\pi}{a_2} \right)^{d/2}\frac{a_1}{\lambda}< \infty.
   \eeqnn

Similarly for $\partial_{x_i}Q^{\lambda}(x)$, by first using (\ref{eqn:interchange})
and then setting instead $m=1$, $y=0$, $r=0$ and $s=1$ in (\ref{LSU}):
    \beqnn
   \int_{\R^d}|\partial_{x_i}Q^{\lambda}(x)| dx
   \leq \int_0^{\infty}e^{-\lambda t} \frac{a_1}{\sqrt{t}}
   \prod_{i=1}^{d} \left(\int_{\R} \frac{1}{\sqrt{t}}
   \exp\{ - a_2 \frac{x_i^2}{t}\}dx_i \right)dt
   = \left( \frac{\pi}{a_2}
   \right)^{d/2}\frac{a_1\sqrt{\pi}}{\sqrt{\lambda}} < \infty.
  \eeqnn

  In the case $d=1$,  $\partial_{x}Q^{\lambda}(x)$ is also square-integrable: the use of
  (\ref{eqn:interchange}) and (\ref{LSU}) yields
  \beqnn
   \int_{\R}|\partial_{x}Q^{\lambda}(x)|^2 dx
   & \leq & \int_{\R} \left | \int_0^{\infty}e^{-\lambda t}
   \frac{a_1}{t} \exp{\{- a_2 \frac{x^2}{t}\}}dt
  \right |^2 dx  \\
   & = &  \int_0^{\infty} \int_0^{\infty}e^{-\lambda s}e^{-\lambda t} \frac{a_1^2}{st}
  \sqrt{ \frac{\pi st}{a_2(s+t)} } ds dt \\
   & \leq &  a_1^2  \sqrt{ \frac{\pi}{2a_2} } \int_0^{\infty} \int_0^{\infty}e^{-\lambda s}e^{-\lambda t} s^{-3/4} t^{-3/4} ds dt
     = a_1^2 \Gamma^2(1/4) \sqrt{ \frac{\pi}{2a_2\lambda} } < \infty,
  \eeqnn
  by first expanding the square and changing the order of integration,
  then using elementary inequality $2\sqrt{ st } \le s+t$ and the classical Gamma function $\Gamma(x)$.

    Finally $Q^{\lambda}(x)$ is square-integrable when $d\le3$: set $y=0$, $r=0$, $s=0$ in (\ref{LSU}) to get
  \beqnn
   \int_{\R^d}|Q^{\lambda}(x)|^2 dx
   & \leq &  \int_{\R^d} \left | \int_0^{\infty}e^{-\lambda t}
   \frac{a_1}{t^{d/2}} \exp{\{- a_2 \frac{|x|^2}{t}\}}dt
   \right |^2 dx \\
   & = &  \int_0^{\infty} \int_0^{\infty}e^{-\lambda s}e^{-\lambda t} \frac{a_1^2}{s^{d/2} t^{d/2}}
  \left( \frac{\pi st}{a_2(s+t)} \right)^{d/2} ds dt \\
   & \leq &  a_1^2   \left( \frac{\pi}{2a_2} \right)^{d/2}
    \int_0^{\infty} \int_0^{\infty}e^{-\lambda s}e^{-\lambda t} s^{-d/4} t^{-d/4} ds dt \\
   & = &  \frac{a_1^2}{\lambda^2} \left( \frac{\pi\lambda}{2a_2} \right)^{d/2}  \Gamma^2(1-d/4)< \infty ,
  \eeqnn
 by the same argumentation, provided $d\le3$ for the finiteness. 
 
For parts (iv) and (v), the proofs will benefit from a bit of simplification in the notation.  
By Lemma \ref{FellerProp}, there is a constant $k>0$ such that 
   \beqlb\label{eqn:interchange2}
    |\partial_{x_i}Q^{\lambda}(x)|
  &=&   |\partial_{x_i}\int_0^{\infty}e^{-\lambda t}q_t^{1}(0,x)dt|  
       =|\int_0^{\infty}e^{-\lambda t}\partial_{x_i}q_t^{1}(0,x)dt|     \nonumber \\
  &\leq& k \int_0^{\infty}e^{-\lambda t} \frac{1}{\sqrt{t}} \varphi_{t/2a_2}(x)dt 
      = k' \int_0^{\infty}e^{-\lambda' t} \frac{1}{\sqrt{t}} \varphi_{t}(x)dt <\infty
 \eeqlb
 holds for all $x \in \R^d \smallsetminus \{0\}$ and any $i \in \{1,2,\cdots,d\}$, 
 where $\varphi_t$ is the transition density of the standard $d$-dimensional Brownian motion,  
 $k'=k\sqrt{2a_2}$ and $\lambda'=2a_2\lambda$. Since only finiteness is sought for every $\lambda>0$, 
 we drop the prime by setting $a_2=1/2$. 

Let us consider (iv). Observe, with a different constant $k_1>0$, 
$$
 \int_{\R^d}| \partial_{x_i}Q^{\lambda}(x-y) |\mu_0(dy)  
 \leq k_T + k \sup_{x \in \R^d} \sup_{0 < t \leq T}\< \varphi_t(x-y), \mu_0(dy)\>\int_0^{T}  e^{- \lambda t} \frac{1}{\sqrt{t}} dt,
 $$
  where
  \beqlb \lab{PL2}
 k_T & := &  k \int_{T}^{\infty}  e^{- \lambda t} \frac{1}{\sqrt{t}}\< \varphi_t(x-y), \mu_0(dy) \>dt         \nonumber \\
         & \leq &      k \int_{T}^{\infty}  e^{- \lambda t} \frac{1}{\sqrt{t}} \frac{1}{(2 \pi T)^{d/2}}\<\exp\{ - \frac{|x-y|^2}{2t}\}, \mu_0(dy) \>dt         \nonumber \\                                                               
         & \leq &\frac{k}{(2 \pi T)^{d/2}} \sup_{x \in \R^d}\<I_a(x-y), \mu_0(dy)\>[ \sup_{x-y \in \R^d} \exp\{ - \frac{|x-y|^2}{2T}\}I^{-1}_a(x-y)] \int_T^{\infty}\frac{e^{- \lambda t}}{\sqrt{t}}dt \nonumber \\
         & < & \infty.
   \eeqlb
   Thus, $ \partial_{x_i}Q^{\lambda}(x-\cdot) \in L^1(\mu_0)$ holds for every $\lambda>0$, 
   using Hypotheses $\ref{hyp:basicassumpGauss}$ and $\ref{hyp:basicassumpUniformInteg}$. \\
   
   Next, we consider (v). Using inequality (\ref{LSU}), observe
   \beqlb \lab{QL2}
 \int_{\R^d}[Q^{\lambda}(x-y)]^2 \mu_0(dy)  & = & \< \int_{0}^{\infty}  e^{- \lambda u} q_u^{1}(x-y)du \int_{0}^{\infty}  e^{- \lambda v} q_v^{1}(x-y)dv, \mu_0(dy) \>       \nonumber \\
         & = &      \int_{0}^{\infty}   \int_{0}^{\infty}   e^{- \lambda (u + v)} \< q_u^{1}(x-y)q_v^{1}(x-y), \mu_0(dy)\>  du dv       \nonumber \\
           & \leq &   k_1   \int_{0}^{\infty}   \int_{0}^{\infty}   e^{- \lambda (u + v)} \<  \varphi_u(x-y) \varphi_v(x-y), \mu_0(dy)\>  du dv       \nonumber \\
         & \leq &   k_1   \int_{0}^{\infty}   \int_{0}^{\infty}   e^{- \lambda (u + v)}    \<\frac{1}{(2 \pi)^{d/2}} \frac{1}{(u+v)^{d/2}}\<  \varphi_{\tau}(x-y),     \mu_0(dy)\>  du dv       \nonumber \\
         &  \leq &  k_1(T_{(i)} + T_{(ii)} + T_{(iii)} + T_{(iv)}),
         \eeqlb
   where  $\tau := \frac{uv}{u+v}$ and
   \beqlb \lab{Ti}
   T_{(i)} := \int_{0}^{T}   \int_{0}^{T}   e^{- \lambda (u + v)}    \<\frac{1}{(2 \pi)^{d/2}} \frac{1}{(u+v)^{d/2}}\<  \varphi_{\tau}(x-y),     \mu_0(dy)\>  du dv,
   \eeqlb
    \beqlb \lab{Ti}
   T_{(ii)} := \int_{T}^{\infty}   \int_{0}^{T}   e^{- \lambda (u + v)}    \<\frac{1}{(2 \pi)^{d/2}} \frac{1}{(u+v)^{d/2}}\<  \varphi_{\tau}(x-y),     \mu_0(dy)\>  du dv,
   \eeqlb
    \beqlb \lab{Ti}
   T_{(iii)} := \int_{0}^{T}   \int_{T}^{\infty}   e^{- \lambda (u + v)}    \<\frac{1}{(2 \pi)^{d/2}} \frac{1}{(u+v)^{d/2}}\<  \varphi_{\tau}(x-y),     \mu_0(dy)\>  du dv,
   \eeqlb
    \beqlb \lab{Ti}
   T_{(iv)} := \int_{T}^{\infty}   \int_{T}^{\infty}   e^{- \lambda (u + v)}    \<\frac{1}{(2 \pi)^{d/2}} \frac{1}{(u+v)^{d/2}}\<  \varphi_{\tau}(x-y),     \mu_0(dy)\>  du dv.
   \eeqlb
 Then, we have
 \beqlb \lab{Ti}
   T_{(i)}  :=  & & \int_{0}^{T}   \int_{0}^{T}   e^{- \lambda (u + v)}    \<\frac{1}{(2 \pi)^{d/2}} \frac{1}{(u+v)^{d/2}}\<  \varphi_{\tau}(x-y),     \mu_0(dy)\>  du dv \nonumber \\
   \leq & & k_2  \int_0^{\pi/2} \int_0^T e^{- \lambda r} \frac{1}{r^{d/2}} 2 r \sin \theta \cos \theta d \theta dr \nonumber \\
   = && k_2 (\sin \theta)|_0^{\pi/2} \int_0^T e^{- \lambda r} r^{1 - d/2}dr = k_2 \int_0^T e^{- \lambda r} r^{1 - d/2}dr< \infty,
   \eeqlb
   if $d=1,2,3$, where
   \[
   k_2 :=  \frac{1}{(2 \pi)^{d/2}} \sup_{x \in \R^d}\sup_{0 < \tau \leq T} \<  \varphi_{\tau}(x-y),     \mu_0(dy) \>.
   \]
   
Similarly one gets successively
 \beqlb \lab{Tii}
   T_{(ii)} & & := \int_{T}^{\infty}   \int_{0}^{T}   e^{- \lambda (u + v)}    \<\frac{1}{(2 \pi)^{d/2}} \frac{1}{(u+v)^{d/2}}\<  \varphi_{\tau}(x-y),     \mu_0(dy)\>  du dv \nonumber \\
   & & \leq \sup_{x \in \R^d} \sup_{0 < \tau \leq T}\<  \varphi_{\tau}(x-y),     \mu_0(dy)\>  \int_{T}^{\infty}   \int_{0}^{T}   e^{- \lambda (u + v)}    \<\frac{1}{(2 \pi)^{d/2}} \frac{1}{(u+v)^{d/2}} du dv \nonumber \\
   & &  < \infty;
   \eeqlb
    \beqlb \lab{Tiii}
   T_{(iii)} & & := \int_{0}^{T}   \int_{T}^{\infty}   e^{- \lambda (u + v)}    \<\frac{1}{(2 \pi)^{d/2}} \frac{1}{(u+v)^{d/2}}\<  \varphi_{\tau}(x-y),     \mu_0(dy)\>  du dv \nonumber \\
   & & \leq \sup_{x \in \R^d} \sup_{0 < \tau \leq T}\<  \varphi_{\tau}(x-y),     \mu_0(dy)\>  \int_{0}^{T}   \int_{T}^{\infty}   e^{- \lambda (u + v)}    \<\frac{1}{(2 \pi)^{d/2}} \frac{1}{(u+v)^{d/2}} du dv \nonumber \\
   & &  < \infty;
   \eeqlb
    \beqlb \lab{Tiv}
   T_{(iv)} & & := \int_{T}^{\infty}   \int_{T}^{\infty}   e^{- \lambda (u + v)}    \<\frac{1}{(2 \pi)^{d/2}} \frac{1}{(u+v)^{d/2}}\<  \varphi_{\tau}(x-y),     \mu_0(dy)\>  du dv \nonumber \\
   & & \leq \sup_{x \in \R^d} \sup_{0 < \tau \leq T}\<  \varphi_{\tau}(x-y),     \mu_0(dy)\>  \int_{T}^{\infty}   \int_{T}^{\infty}   e^{- \lambda (u + v)}    \<\frac{1}{(2 \pi)^{d/2}} \frac{1}{(u+v)^{d/2}} du dv \nonumber \\
   & &  < \infty.
   \eeqlb
So this proves that $Q^{\lambda} \in L^2(\mu_0 )$ for $\mu_0 \in M_a(\R^d)$, $a > 0$ and  $d=1,2,3$.

\subsection{Proof of upper bound (\ref{DI13})}\label{app:pf_DI13}
 \proof
 By the monotonicity in $d$ of the integrands in (\ref{DI13}), when $r$, $u$ and $v$ are arbitrarily fixed,
 the bound $\chi_{2}(s)\le \chi_{1}(s)+\chi_{3}(s)$ is valid for all values of $s\ge0$. (By the way, a simple integration by parts
 also yields $\chi_{4}(s)=\infty$ everywhere.)  Next, because of inequalities
 \beqlb
  \int_0^{s}  \frac{1}{(\sqrt{2r +u + v})^{3}} dr
  = -  \frac{1}{\sqrt{2s +u +v}} + \frac{1}{\sqrt{u +v}}
  \le \frac{1}{\sqrt{u +v}} \nonumber
 \eeqlb
 and
  \beqlb
  \int_0^{s}  \frac{1}{\sqrt{2r +u + v}} dr
  \le \frac{s}{\sqrt{u +v}}, \nonumber
 \eeqlb
 we can also see that the integrals $\chi_{1}(s)$,  $\chi_{2}(s)$ and $\chi_{3}(s)$ are each bounded by $(s+1)$ times
  \beqlb
 \int_0^{\infty} \int_0^{\infty} e^{- \lambda (u+v)} \frac{1}{\sqrt{uv}} \frac{1}{\sqrt{u + v}} du dv. \nonumber
 \eeqlb

 By first using the transformation $u= w^2$ and $v=z^2$, followed by a polar coordinate transformation,
 this last integral becomes successively
 \beqlb
 \int_0^{\infty} \int_0^{\infty} e^{- \lambda (w^2 + z^2)} \frac{4}{\sqrt{w^2 + z^2}} dw dz
 = 4 \int_0^{\frac{\pi}{2}}  \int_0^{\infty} e^{- \lambda r^2} dr d\theta
  = \pi \sqrt{\frac{\pi}{\lambda}}, \nonumber
 \eeqlb
 yielding the common upper bound (\ref{DI13}) for $\chi_{d}(s)$, valid for each of $d=1$, 2 and 3.
\qed

\subsection{Proof of upper bound (\ref{DI1a})}\label{app:pf_DI1a}
 \proof 
 For any $x,y\in\R^d$ denote by $\Psi_{xy}^{\epsilon,\delta}\in C_b({(\R^d)^2})\cap L^1(\R^{2d})\cap L^1(\mu_0^2)$ the function
\[
\Psi_{xy}^{\epsilon,\delta}(w_1,w_2):=\Phi_{xy}^{\epsilon,\delta}(w_1)\times\Phi_{xy}^{\epsilon,\delta}(w_2)
\] 
with
\[
\Phi_{xy}^{\epsilon,\delta}(w):=h_p(y-w) [\partial_pQ^{\lambda}_{\epsilon}(x-w) - \partial_pQ^{\lambda}_{\delta}(x-w)].
\] 
Using (\ref{LSU2}), the coarse upper bound
\[
|\partial_pQ^{\lambda}_{\epsilon}(x) - \partial_pQ^{\lambda}_{\delta}(x) | 
\leq c(e^{\lambda \epsilon} -1) \int_{\epsilon}^{\infty}e^{- \lambda u} \frac{1}{\sqrt{u}}g_{u}(x)du  
+ ce^{\lambda\epsilon}\int_{0}^{\epsilon}e^{- \lambda u} \frac{1}{\sqrt{u}}g_{u}(x)du,
\]
valid for any $0 < \delta < \epsilon < \infty$, combined with (\ref{DI10}) yields a constant $c > 0$ such that 
\beqnn
  \int_{\r^d} | \Psi_{xy}^{\epsilon,\delta}(w_1,w_2) | dy 
   & \le & c(e^{\lambda \epsilon} -1)^2 \int_{\epsilon}^{\infty}\int_{\epsilon}^{\infty} 
   e^{- \lambda (u+v)} \frac{1}{\sqrt{uv}} g_u(x-w_1)g_v(x-w_2) dudv \\
  & & + 2ce^{\lambda \epsilon}(e^{\lambda \epsilon} -1) \int_0^{\epsilon}\int_{\epsilon}^{\infty} 
   e^{- \lambda (u+v)} \frac{1}{\sqrt{uv}} g_u(x-w_1)g_v(x-w_2) dudv \\
  & & + c(e^{\lambda \epsilon})^2 \int_0^{\epsilon}\int_0^{\epsilon}
   e^{- \lambda (u+v)} \frac{1}{\sqrt{uv}} g_u(x-w_1)g_v(x-w_2) dudv 
 \eeqnn
 for all $0 < \delta < \epsilon < \infty$, with the four terms merged into three. 
 Proceeding as in the argument leading to (\ref{DI12}), we get a new constant $c > 0$,
 independent of $\lambda$, such that
\beqnn
  \int_0^t ds \int_{\r^d} dy \int_0^s \< P^1_{s-r}\Phi_{12}^2P^2_r\Psi_{xy}^{\epsilon,\delta}, \mu_0 \> dr  
  & \leq & c (e^{\lambda \epsilon} -1)^2 \chi_{d}^{\infty,\infty}(t) \\
  & & + 2c e^{\lambda \epsilon}(e^{\lambda \epsilon} -1) \chi_{d}^{\epsilon,\infty}(t) 
         + c e^{2\lambda \epsilon} \chi_{d}^{\epsilon,\epsilon}(t) 
 \eeqnn 
with $\chi_{d}^{\epsilon,\eta}$ given by
\beqlb \lab{chiepsdelta}
 \chi_{d}^{\epsilon,\eta}(s) \equiv  \int_0^{\epsilon} \int_0^{\eta} e^{- \lambda (u+v)} \frac{1}{\sqrt{uv}}
  \int_0^s \frac{1}{(\sqrt{2r +u + v})^{d}} dr du dv
\eeqlb
and going to $0$ with either $\epsilon$ or $\eta$ (or both) when $d\le3$, 
using the finiteness in (\ref{DI13}); and
 \beqnn
 \int_0^t ds \int_{\r^d} dy \< P^2_s\Psi_{xy}^{\epsilon,\delta}, \mu_0 \>  
 & \leq & c(e^{\lambda \epsilon} -1)^2  \int_0^t ds 
 \int_{\epsilon}^{\infty}\int_{\epsilon}^{\infty} e^{- \lambda (u+v)} \frac{1}{\sqrt{uv}} dudv \\
 & & + 2ce^{\lambda \epsilon}(e^{\lambda \epsilon} -1)  \int_0^t ds 
  \int_0^{\epsilon}\int_{\epsilon}^{\infty}  e^{- \lambda (u+v)} \frac{1}{\sqrt{uv}} dudv \\ 
 & & + ce^{2\lambda \epsilon}  \int_0^t ds 
 \int_0^{\epsilon}\int_0^{\epsilon} e^{- \lambda (u+v)} \frac{1}{\sqrt{uv}} dudv,
 \eeqnn
which goes to $0$ with $\epsilon$, at a rate of $\epsilon$ (because of the third and asymptotically largest term), 
since $\int_{0}^{\infty} e^{- \lambda u} \frac{1}{\sqrt{u}} du < \infty$ and 
$e^{\lambda \epsilon} -1\le (e^{\lambda} -1)\epsilon$ when $0\le\epsilon\le1$.  
The vanishing rate for the other term follows from the calculations in Subsection \ref{app:pf_DI13},
which yield
\beqlb 
  \left| \chi_{d}^{\epsilon,\eta}(s) \right|
  \le 2\pi (s+1) \int_0^{\sqrt{\epsilon^2+\eta^2}} e^{- \lambda r^2} dr   
  \le 2\pi (s+1) \sqrt{\epsilon^2+\eta^2}, \nonumber 
 \eeqlb
hence $\chi_{d}^{\epsilon,\epsilon}(s)$ goes to $0$ at a rate of $\epsilon$. 
\qed

\noindent {\bf Acknowledgements.}
The third author would also like to thank Professor Daniel Dugger for assistance in  
accessing the University of Oregon research resources. 


%
%






\end{document}